\documentclass[11pt, a4paper]{article}

\usepackage[utf8]{inputenc}
\usepackage[T1]{fontenc}
\usepackage{lmodern}

\usepackage[margin=1in]{geometry}
\usepackage{setspace}
\usepackage{titlesec}
\titleformat*{\section}{\large\bfseries\boldmath}
\titleformat*{\subsection}{\normalsize\bfseries\boldmath}
\titleformat*{\paragraph}{\normalsize\bfseries\boldmath}
\usepackage{amsmath, amssymb, amsthm, amsfonts, mathtools, xspace}
\usepackage{booktabs,microtype,needspace,placeins}
\allowdisplaybreaks

\usepackage{tcolorbox}
\tcbuselibrary{breakable}
\usepackage[noEnd=true,indLines=true]{algpseudocodex}
\tikzset{algpxIndentLine/.style={draw=black,line width=0.4pt}}
\newtcolorbox[auto counter]{boxedalgorithm}[2][]{%
  breakable,
  colback=white,
  colframe=black,
  colbacktitle=white,
  coltitle=black,
  fonttitle=\bfseries,
  sharp corners,
  boxrule=0.4pt,
  before skip=\medskipamount,
  after skip=\medskipamount,
  title={Algorithm~\thetcbcounter. #2},
  #1
}
\usepackage{graphicx}
\usepackage{tikz}
\usetikzlibrary{patterns}
\definecolor{elsBlue}{HTML}{0072B2}
\definecolor{elsOrange}{HTML}{D55E00}
\definecolor{elsGray}{HTML}{E8ECEF}
\definecolor{elsInk}{HTML}{252C31}
\definecolor{elsMid}{HTML}{697782}
\definecolor{elsGuide}{HTML}{D9E0E5}
\tikzset{
  els/base/.style={
    font=\small,
    text=elsInk,
    line cap=round,
    line join=round
  },
  els/axis/.style={draw=elsMid,line width=.45pt},
  els/guide/.style={draw=elsGuide,line width=.35pt},
  els/curve/.style={line width=1.1pt},
  els/block/.style={line width=.5pt},
  els/leader/.style={line width=.45pt}
}
\usetikzlibrary{arrows.meta}
\tikzset{
  els/mapping/.style={
    -{Stealth[length=3pt]},
    draw=elsMid!65,
    line width=.45pt
  }
}
\usepackage[colorinlistoftodos, textsize=footnotesize]{todonotes}
\usepackage[authoryear,round]{natbib}

\definecolor{darkblue}{rgb}{0.0, 0.0, 0.55}
\usepackage[pagebackref=true, colorlinks=true, allcolors=darkblue]{hyperref}
\renewcommand*{\backref}[1]{}
\renewcommand*{\backrefalt}[4]{%
    \ifcase #1%
    \or #2%
    \else #2%
    \fi%
}

\makeatletter
\newtheoremstyle{smallcaps}
  {}
  {}
  {\itshape}
  {}
  {\scshape}
  {.}
  { }
  {}

\renewenvironment{proof}[1][\proofname]{\par
  \pushQED{\qed}%
  \normalfont \topsep6\p@\@plus6\p@\relax
  \trivlist
  \item[\hskip\labelsep
        \scshape
    #1\@addpunct{.}]\ignorespaces
}{%
  \popQED\endtrivlist\@endpefalse
}
\makeatother
\theoremstyle{smallcaps}
\newtheorem{theorem}{Theorem}[section]
\newtheorem{lemma}[theorem]{Lemma}

\newtheorem{claim}[theorem]{Claim}

\newcommand{\eps}{\epsilon}
\newcommand{\ind}{\text{\usefont{U}{dsrom}{m}{n}1}}
\newcommand{\OPT}{\operatorname{OPT}}
\newcommand{\OPTcyc}{\operatorname{OPT}_{\mathrm{cyc}}}
\newcommand{\OPTdyn}{\operatorname{OPT}_{\mathrm{dyn}}}
\newcommand{\dd}{\mathrm{d}}
\newcommand{\bbN}{\mathbb{N}}

\newcommand{\mysetup}{\mathrm{setup}}
\newcommand{\myhold}{\mathrm{hold}}
\newcommand{\mylast}{\mathrm{last}}

\title{\boldmath\textbf{A Constant-Factor Approximation for \\
Economic Lot Scheduling}}
\author{Danny Segev\thanks{School of Mathematical Sciences, Tel Aviv University,
Tel Aviv 69978, Israel. Email: \texttt{segevdanny@tauex.tau.ac.il}. Research supported by Israel Science Foundation grant 767/26.}}
\date{}

\begin{document}

\maketitle
\thispagestyle{empty}

\begin{abstract}
\noindent Economic lot scheduling is one of the most fundamental and extensively studied problems in production scheduling and inventory management.
This setting asks how multiple commodities sharing a common machine should be
produced to meet continuous demand while minimizing the long-run average cost of setups and inventory holding. Since the 1950s, it has motivated hundreds of studies on production frequencies, cyclic
schedules, feasibility conditions, and various optimization methods. Despite this sustained attention, provably good approximation guarantees for the general problem formulation,
measured against unrestricted dynamic policies, have remained elusive.

\medskip \noindent The main contribution of this paper resides in proposing the first constant-factor approximation for the economic lot scheduling problem. This performance guarantee holds in relation to the unrestricted dynamic optimum, without relying on any parametric assumptions beyond those ensuring feasibility. Our main technical innovation is a new convex relaxation that captures the inventory consequences of competition for machine time, combined with a priority-based scheduling rule
whose analysis relates machine occupation to these inventory requirements. Concurrently, we overcome the potentially exponential size of non-trivial cyclic schedules, showing that the resulting policy admits a polynomial-size description that allows its
individual production decisions to be efficiently recovered.
\end{abstract}

\vspace{1em}

\newpage
\pagestyle{empty}
\tableofcontents
\newpage


\pagestyle{plain}
\pagenumbering{arabic}
\section{Introduction} \label{sec:intro}

Economic lot scheduling is one of the most fundamental problems
at the intersection of production scheduling and inventory management.
In this setting, we wish to coordinate the production of several commodities on a common machine, where each commodity faces a constant demand rate and
can be produced at a faster rate. However, at any point in time, only one commodity may be processed.
Starting a production run incurs a setup cost and may require setup
time, while inventory carried between runs incurs holding costs. Our objective
is to meet demand throughout a continuous-time infinite planning horizon while minimizing the long-run
average cost of setups and inventory holding. For a single commodity, the inherent tradeoff we are dealing with is well understood: Longer
runs reduce the frequency of setups but require more inventory between
successive runs. However, with multiple commodities sharing a single machine, each run
also prevents other commodities from being produced, and their inventories
must cover demand throughout this interruption, meaning that economical production
frequencies cannot be independently chosen. 

Economic lot scheduling has attracted
sustained attention since the 1950s, with historical landmarks including
the joint treatment
of lot sizing and schedule feasibility by \citet{Rogers1958}, the early work of \citet{Maxwell1964} on cyclic production schedules, and the dynamic programming
approach of \citet{Bomberger1966}. Subsequent research has developed a
substantial body of work on production frequencies, cyclic policies,
feasibility conditions, and various computational methods. The scale of this literature is illustrated by the 242 papers
analyzed by \citet{BeckGlock2020}, the 141 common-cycle studies
reviewed by \citet{NobilEtAl2026}, and our most recent Google Scholar
search (November 2026), which returned approximately 341 results with the phrase
``economic lot scheduling'' in their titles. Reviews of this
literature include \citet{Elmaghraby1978}, \citet{Carstensen1999},
\citet{ChanEtAl2013}, \citet{SantanderMercadoJubizDiaz2016}, and
\citet{BeckGlock2020}. Detailed textbook treatments are available in
\citet[Sec.~5.6]{Zipkin2000}, \citet[Chap.~7]{Pinedo2005},
\citet[Sec.~7.2.1]{Axsater2015},
\citet[Sec.~4.9]{NahmiasOlsen2015}, and
\citet[Sec.~10.6.1]{SilverPykeThomas2017}.

\paragraph{What existing approaches accomplish.}
Many of these research efforts obtain tractable solutions by imposing
structural restrictions on the production schedule. Two extensively studied classes are common-cycle (or rotation) policies, where 
each commodity is produced once per period, and basic-period policies, allowing
different production intervals that are integer multiples of a
common time unit. Well-known contributions along these lines include the dynamic
programming approach of \citet{Bomberger1966} and the basic-period heuristic of \citet{DollWhybark1973} for choosing production
frequencies and a basic cycle time. Selecting frequencies
alone, however, does not ensure that the corresponding runs can
be scheduled without overlap. For rotation cycles,
\citet{JonesInman1989} derived performance bounds and identified
parameter regimes in which producing each commodity once per
cycle is optimal or nearly optimal. More flexible cyclic policies allow each commodity's lot sizes
and the intervals between its successive runs to vary within
the cycle. \citet{Dobson1987} developed a formulation that
exploits this flexibility to construct feasible schedules while
explicitly accounting for setup times. He derived feasibility
conditions and quantified the effect of small schedule
adjustments on cost. For a
prescribed sequence in this formulation, \citet{Zipkin1991}
optimized production durations, idle times, and cycle length via
parametric quadratic programming. This approach separates the continuous
timing decisions from the combinatorial choice of a production
sequence. Regularity in production intervals also plays a role in
approximation methods. For example, \citet{Roundy1989} developed
power-of-$2$ rounding methods for single-resource lot sizing, which can be viewed as a relaxation of economic lot scheduling. However, these ideas do not provide a corresponding approximation for a feasible production schedule. \citet{AdelmanBarz2014} derived tractable lower bounds through
an approximate dynamic programming formulation and showed how
this framework recovers earlier lower-bound formulations. In terms of rigorous complexity results, NP-hardness arises even in restricted cyclic
formulations, as shown by \citet{Hsu1983}.
More generally, \citet{GallegoShaw1997} established strong
NP-hardness of optimizing over general cyclic schedules,
where the production sequence, lot sizes, and intervals
between successive runs are all decision variables.

\paragraph{Motivating questions.} These contributions provide
exact methods, practical heuristics, and performance guarantees
under specific structural restrictions or parametric assumptions. However, to date, we are unware of any provably good polynomial-time approximation guarantees against unrestricted dynamic policies for the general problem formulation, where commodities have arbitrary production rates, setup costs, setup times, and holding costs. Consequently, our work is motivated by three long-standing open questions:
\begin{itemize}
    \item \emph{Approximation guarantees.}
    Can the minimum long-run average cost be efficiently approximated within
    any non-trivial  factor? The benchmark we are referring to ranges
    over all feasible dynamic policies, without requiring a fixed
    production sequence, uniform lot sizes, or periodic repetition.

    \item \emph{Capturing competition for machine time.}
    Can tractable relaxations account for the inventory that must be carried by each
    commodity  while the machine sets up for or produces
    other commodities? Can these interactions be
    captured by a mathematical program that yields provably good lower bounds on the dynamic optimum?

    \item \emph{Policy construction and representation.}
    Can such a program guide the construction of a feasible
    policy with comparable cost? Even when a cyclic policy suffices,
    its period may contain too many production runs to list
    explicitly. Can we instead compute a compact description
    from which the policy's activity at any given time can be efficiently
    determined?
\end{itemize}

\paragraph{Informal contributions.}
We prove that the economic lot scheduling problem admits a constant-factor
approximation that can be efficiently computed and compactly represented. This performance guarantee is measured against the unrestricted dynamic optimum,
allowing arbitrary production sequences and run durations. Specifically, we develop a deterministic polynomial-time
approach for computing a lower bound on this optimum together
with a feasible cyclic policy whose cost is at most $18$ times
that bound. As shown in the sequel, for certain instances, any non-trivial cyclic policy requires
exponentially many production runs per period. We overcome this
obstacle by constructing a polynomial-size description from which
the policy's activity at any time, including its starting
and ending points, can be efficiently recovered. The resulting
policy can therefore be executed without ever expanding its full cycle. At the heart of our approach are a new convex relaxation that accounts for
the inventory required when commodities compete for machine time,
a priority-based rule that converts its solution into a feasible
production schedule, and an implicit representation that supports
efficient access to this schedule.

\paragraph{Outline.} Moving forward, Section~\ref{sec:model} provides a formal mathematical description of the economic lot scheduling model,
feasible dynamic policies, and their long-run cost objective.
Section~\ref{sec:cyclic-queries} introduces cyclic policies,
explains why a compact representation is necessary, and defines
the notion of activity queries.
Section~\ref{sec:contributions} states the main contributions of this paper and outlines the central ideas behind
our convex relaxation, policy construction, and activity queries.

\subsection{Model formulation}\label{sec:model}

In the economic lot scheduling problem, we are given $n$ commodities, indexed by $1,\ldots,n$, which share a single production machine over the continuous-time planning horizon $[0,\infty)$. Each commodity $i$ is associated with a stationary demand rate $d_i>0$, a production rate $p_i>d_i$, a setup cost $K_i>0$, a setup time $s_i\ge0$, and a holding cost $h_i>0$ per unit of inventory per unit of time. The precise roles of these input parameters will be explained below, noting that for computational purposes, they are all assumed to be rationally valued.

\paragraph{Capacity requirements.} Beyond these elementary sign and rationality conditions, we impose only an
easily verifiable feasibility condition. Specifically, let
$\rho_i=\frac{d_i}{p_i}$ be the fraction of machine time required
to meet the demand for commodity $i$, and let
$R=\sum_{i\in[n]}\rho_i$ be the aggregate fraction of machine
time required for production. As we prove in Appendix~\ref{app:feasibility}, a feasible
policy exists if and only if either $R<1$, or $R=1$ and all setup times
are zero. For ease of presentation,
our analysis considers $R<1$, leaving positive capacity
for setups and idling, without requiring the remaining capacity
$\kappa=1-R$ to be lower-bounded by any fixed positive constant. With minor adaptations along the way, our approach extends to the case of $R=1$ and $s_1=\cdots=s_n=0$; we omit these details to avoid redundancy.

\paragraph{Production policies.} A dynamic production policy $\pi$ specifies a finite initial
inventory vector $I^\pi(0)\in\mathbb R_+^n$ and an
infinite sequence of production runs. Each run is represented by a triple $(i_k,a_k,b_k)$, where $i_k$ is the commodity
produced, $a_k$ is the time at which its setup begins, and $b_k$ is the
time at which its production ends. During this run, the machine is set up for
commodity $i_k$ on $[a_k,a_k+s_{i_k})$, incurring the setup cost $K_{i_k}$, and subsequently produces  $i_k$ on $[a_k+s_{i_k},b_k)$ at rate
$p_{i_k}$. We require $a_k\ge 0$ and
$a_k+s_{i_k}<b_k \le a_{k+1}$ for every $k$, so that each run has positive production
time, and neither its setup nor its production can overlap the succeeding run.
At all other times, the machine is idle. We enforce that $(a_k)_{k \in \bbN}$ tends to infinity, meaning that only finitely many
setups begin in any bounded interval.

\paragraph{Inventory dynamics and feasibility.} For a given policy $\pi$, let $I_i^{ \pi }(t)$ denote the inventory of commodity
$i$ at time $t$, and let
$\mathcal{P}_i = \bigcup_{k:i_k=i} [a_k+s_i,b_k)$ be the set of time intervals at which we produce this commodity.
Production enters inventory continuously at rate $p_i$ on
$\mathcal{P}_i$, whereas demand occurs at rate $d_i$ at all times.
Thus, for the initial stock $I_i^{ \pi }(0)$, our inventory level at any time $t\ge 0$ can be written as
\begin{equation} \label{eqn:inventory_time}
I_i^{ \pi }(t)
~~=~~ I_i^{ \pi }(0) + p_i \cdot |\mathcal{P}_i\cap[0,t)| - d_i t,    
\end{equation}
where $|\mathcal{P}_i\cap[0,t)|$ is the total time spent
producing commodity $i$ before $t$. In particular, $I_i^\pi$ is absolutely
continuous and satisfies $\frac{ \dd I_i^{ \pi }(t) }{ \dd t } = p_i \cdot \ind \{ t \in \mathcal{P}_i \} - d_i$ for almost every $t\ge 0$. Since demand must be met as it occurs, we say that $\pi$ is a feasible policy if $I_i^{ \pi }(t)\ge 0$ for every commodity $i \in [n]$ and every $t\ge 0$.

\paragraph{Cost structure and objective.}
For every commodity $i\in[n]$ and bounded planning horizon $[0,T)$, let $N_i^{\pi}(T)$ be the number of commodity-$i$ runs whose setup begins in $[0,T)$, namely, $N_i^{\pi}(T)
= |\{k : i_k=i, 0\le a_k<T\}|$. Each such setup incurs a cost of $K_i$, charged at its starting time, meaning that our total
setup cost over $[0,T)$ is $\sum_{i\in[n]}K_iN_i^\pi(T)$. Concurrently, at time $t$, the policy $\pi$ holds $I_i^\pi(t)$ units of commodity $i$,
each incurring a holding cost of $h_i$ per unit of time. Thus, the
holding cost rate at time $t$ is $\sum_{i\in[n]}h_iI_i^\pi(t)$, and integrating this rate over $[0,T)$ leads to our total holding cost,
$\int_0^T\sum_{i\in[n]}h_iI_i^\pi(t) \dd t$. By combining these two components, the long-run average cost of a feasible
policy $\pi$ is defined as 
\[ \mathcal C(\pi)
~~=~~ \limsup_{T\to\infty} \frac{1}{T} \cdot 
\left(
  \sum_{i\in[n]} K_i N_i^{\pi}(T)
  + \int_0^T \sum_{i\in[n]} h_i I_i^{\pi}(t) \dd t
\right). \]
Letting $\Pi$ be the class of all feasible dynamic policies, we make use of $\OPTdyn = \inf_{\pi\in\Pi} \mathcal C(\pi)$ to denote the infimum long-run average cost over all such policies.

\subsection{Cyclic policies and activity queries} \label{sec:cyclic-queries}

Cyclic policies will play an important role in our approach, and we therefore proceed by laying down their basics, along with demonstrating the necessity of a compact representation.

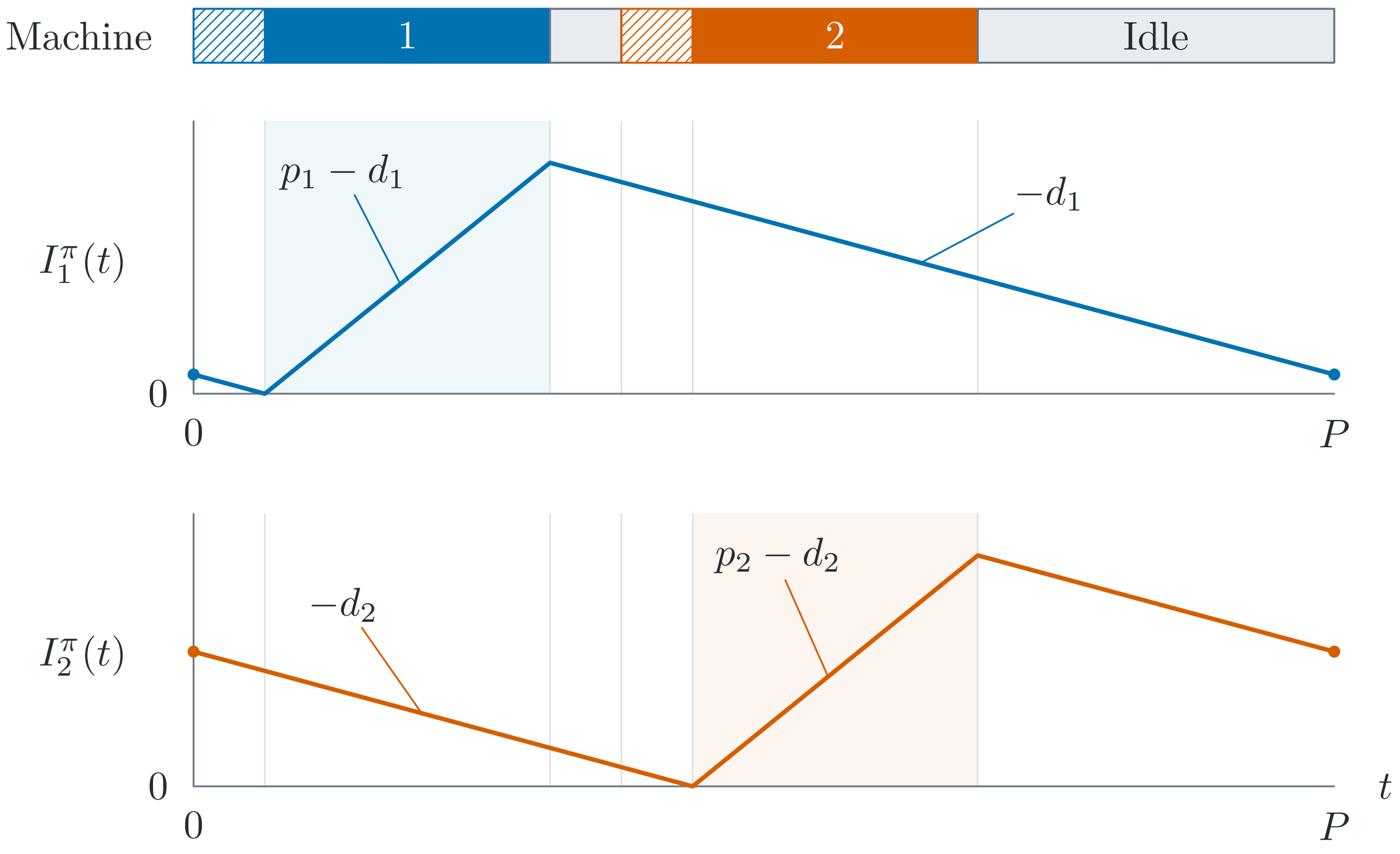
\begin{figure}[htbp!]
\centering
\resizebox{0.95\linewidth}{!}{%
\begin{tikzpicture}[els/base,x=1.25cm,y=1.35cm]


  \path[
    els/block,
    draw=elsBlue,
    pattern=north east lines,
    pattern color=elsBlue
  ] (0,4.70) rectangle (.5,5.05);

  \path[els/block,draw=elsBlue,fill=elsBlue]
    (.5,4.70) rectangle (2.5,5.05);

  \path[els/block,draw=elsMid,fill=elsGray]
    (2.5,4.70) rectangle (3,5.05);

  \path[
    els/block,
    draw=elsOrange,
    pattern=north east lines,
    pattern color=elsOrange
  ] (3,4.70) rectangle (3.5,5.05);

  \path[els/block,draw=elsOrange,fill=elsOrange]
    (3.5,4.70) rectangle (5.5,5.05);

  \path[els/block,draw=elsMid,fill=elsGray]
    (5.5,4.70) rectangle (8,5.05);

  \node[anchor=east] at (-.18,4.875) {Machine};
  \node[text=white] at (1.5,4.875) {$1$};
  \node[text=white] at (4.5,4.875) {$2$};
  \node at (6.75,4.875) {Idle};


  \fill[elsBlue!6]
    (.5,2.55) rectangle (2.5,4.32);

  \foreach \x in {.5,2.5,3,3.5,5.5}
    \draw[els/guide] (\x,2.55)--(\x,4.32);

  \draw[els/axis]
    (0,4.32)--(0,2.55)--(8,2.55);

  \draw[els/curve,elsBlue]
    (0,2.675)--(.5,2.55)--(2.5,4.05)--(8,2.675);

  \fill[elsBlue] (0,2.675) circle (1.5pt);
  \fill[elsBlue] (8,2.675) circle (1.5pt);

  \node[anchor=east] at (-.36,3.40) {$I_1^\pi(t)$};
  \node[anchor=east] at (-.07,2.55) {$0$};
  \node[anchor=north] at (0,2.48) {$0$};
  \node[anchor=north] at (8,2.48) {$P$};

  \node at (1.05,3.99) {$p_1-d_1$};
  \draw[els/leader,elsBlue]
    (1.13,3.84)--(1.45,3.2625);

  \node at (6.00,3.85) {$-d_1$};
  \draw[els/leader,elsBlue]
    (5.75,3.72)--(5.10,3.40);


  \fill[elsOrange!6]
    (3.5,0) rectangle (5.5,1.77);

  \foreach \x in {.5,2.5,3,3.5,5.5}
    \draw[els/guide] (\x,0)--(\x,1.77);

  \draw[els/axis]
    (0,1.77)--(0,0)--(8,0);

  \draw[els/curve,elsOrange]
    (0,.875)--(3.5,0)--(5.5,1.5)--(8,.875);

  \fill[elsOrange] (0,.875) circle (1.5pt);
  \fill[elsOrange] (8,.875) circle (1.5pt);

  \node[anchor=east] at (-.36,.85) {$I_2^\pi(t)$};
  \node[anchor=east] at (-.07,0) {$0$};
  \node[anchor=north] at (0,-.07) {$0$};
  \node[anchor=north] at (8,-.07) {$P$};
  \node[anchor=west] at (8.20,0) {$t$};

  \node at (1.05,1.18) {$-d_2$};
  \draw[els/leader,elsOrange]
    (1.18,1.03)--(1.6,.475);

  \node at (4.10,1.50) {$p_2-d_2$};
  \draw[els/leader,elsOrange]
    (4.15,1.34)--(4.45,.7125);

\end{tikzpicture}%
}
\caption{A feasible cycle for two commodities. Hatched intervals
are setups, numbered intervals are production runs, and gray
intervals are idle. Inventory rises at rate $p_i-d_i$ during
production, falls at rate $d_i$ otherwise, and returns to its
initial level after one period.}
\label{fig:production-inventory-dynamics}
\end{figure}

\paragraph{Cyclic policies.} A policy $\pi$ is cyclic with period $P > 0$ if its
schedule repeats a finite list of runs $((i_k,a_k,b_k))_{k \in [m]}$, where each run satisfies $0\le a_k\le a_k+s_{i_k}<b_k\le P$ and
$b_k\le a_{k+1}$. For every integer $r\ge0$, the runs in
$[rP,(r+1)P)$ are precisely $((i_k,rP+a_k,rP+b_k))_{k \in [m]}$. The machine is idle at all other times, and inventories are carried
across period boundaries without being reset. Since every interval
of length $P$ contains the same amount of production time for each
commodity $i$, at any time $t\ge0$ we have
\[
I_i^\pi(t+P)-I_i^\pi(t)
~~=~~ p_i \cdot \sum_{ \substack{k \in [m]: \\ i_k=i}}(b_k-a_k-s_i)-d_iP .
\]
If this difference were negative, repeated periods would eventually
exhaust the initial stock of $i$. If it were positive, each successive period
would carry that much more stock, leading to an
infinite average holding cost. Thus, for any feasible cyclic policy
of finite long-run average cost, each commodity's production quantity per period equals
its demand, and $I_i^\pi(t+P)=I_i^\pi(t)$ for all $t\ge0$. Figure~\ref{fig:production-inventory-dynamics} illustrates one
period of a feasible two-commodity policy and its inventory
trajectories. Inventory decreases during setups and idle
intervals, and producing either commodity requires the other
to meet demand from stock.

\paragraph{The inevitable complexity of good policies.} Surprisingly, even with only two commodities, cyclic policies with polynomial
approximation guarantees can require exponentially many runs per
period. As an example, in Appendix~\ref{sec:cycle-length}, we
construct instances of binary encoding length $\Theta(k)$ such that, for every $\alpha\ge1$, any cyclic policy of long-run
average cost at most $\alpha\cdot\OPTdyn$ necessarily contains at least
$\frac{4^k}{24\alpha^2}$ runs per period. Thus, even when the
permitted approximation factor is polynomial in the input length,
explicitly listing one period can require exponentially many
entries. 

\paragraph{Activity queries.} This observation motivates us to seek a compact policy description from which
individual production decisions can be efficiently recovered. Such a description should allow individual activities to be identified
without listing all runs. For this purpose, given a
rational time $t\ge0$, an activity query returns:
\begin{itemize}
    \item The activity in progress at time $t$ (setup, production, or idling).
    
    \item The commodity involved, if applicable.

    \item The absolute start and end times of this activity.
\end{itemize}
To avoid ambiguity,  all activity intervals are closed on the left and open on the right.

\subsection{Main results and technical ideas}\label{sec:contributions}

The main contribution of this paper resides in proposing the first constant-factor approximation for the economic lot scheduling problem. Our performance guarantees hold in relation to the unrestricted dynamic optimum, without relying on any parametric assumptions beyond those ensuring feasibility. Concurrently, we overcome the potentially exponential size of approximate cyclic schedules, showing that the resulting policy admits a polynomial-size description that allows its
individual production decisions to be efficiently recovered.
Our main technical innovations are a new convex relaxation that accounts for
the inventory required when commodities compete for machine time,
a priority-based rule that converts its solution into a feasible
production schedule, and an implicit representation that supports
efficient access to this schedule.

\paragraph{Main result 1: Value approximation and compact description.} Along these lines, we devise a computable lower
bound on the unrestricted dynamic optimum and construct a
feasible cyclic policy whose long-run average cost is within a constant factor
of this bound. As formally stated below, the resulting approximation factor is independent of the
number of commodities, the relative magnitudes of their
parameters, and the capacity available for setups and idling. For ease of exposition,  we intentionally do not optimize constants, focusing instead on simplicity and accessibility. 

\begin{theorem} \label{thm:main}
For every feasible instance, we can compute in polynomial time a rational number $L>0$
and a polynomial-size description of a feasible cyclic policy
${\pi}$ satisfying
\[ L ~~\le~~ \OPTdyn ~~\le~~ \mathcal C({\pi}) ~~\le~~ 18L . \]
\end{theorem}

\paragraph{Main result 2: Activity queries.}
The above-mentioned description allows us to efficiently identify the activity performed by our approximate policy ${\pi}$ at any rational time $t\ge0$, including its commodity
and absolute endpoints. The information needed to answer such queries is computed only
once, even when ${\pi}$ contains exponentially many
production runs. To formalize this result, let $P$ be the
policy's period length and let $\Delta$ be a computed positive integer such that all activity endpoints are integer multiples of $\frac{ 1 }{ \Delta }$. Then, we prove that both $P$ and
$\Delta$ have polynomial binary encoding lengths, and moreover, that any activity query has polynomial bit complexity in the encoding length of the instance considered.

\begin{theorem} \label{thm:activity-queries}
With polynomial-time preprocessing, every activity
query for our policy ${\pi}$ can be answered in
$O(n^3\log(P\Delta))$ arithmetic operations.
\end{theorem}

\paragraph{The convex relaxation.} In
Section~\ref{sec:lower}, we develop a convex relaxation that forms the starting point of our algorithmic approach. The resulting formulation can be viewed as strengthening the classical relaxation of 
\citet{Bomberger1966}. For making an easy comparison, the latter can be phrased in terms of run frequencies $\{ f_i \}_{i \in [n]}$ as follows:
\begin{equation} \label{eq:bomberger}
L_{\rm B} ~~=~~ \min \left\{ \sum_{i\in[n]} \left( K_i f_i+\frac{h_i d_i(1-\rho_i)}{2f_i} \right) : \sum_{i\in[n]} s_i f_i \le \kappa, f > 0 \right\}.
\end{equation}
This relaxation accounts for each commodity's own setup and
inventory costs and limits the average capacity consumed by
setups. However, its weakness lies in not charging commodities for the
inventory they must carry during other production runs.
Indeed, Appendix~\ref{sec:gap} shows that the ratio $\frac{ \OPTdyn }{ L_{\rm B} }$ is generally
unbounded, even for two-commodity instances. Our formulation introduces the missing interaction by taking into account,
for every ordered pair of commodities $i \neq j$, the inventory $i$ must carry while the machine sets up for or produces $j$.
Inventory balance and the restriction to a single machine
yield convex constraints on these quantities, thereby bounding the inventory that any group of commodities must carry
while the machine serves commodities outside that group. In this context, a cyclic-approximation argument  ensures that the
resulting lower bound applies to arbitrary dynamic policies,
even though this analysis is derived by examining cyclic policies.

\paragraph{Constructing and analyzing a cyclic policy.} In
Section~\ref{sec:construction}, we convert an approximate solution of our convex relaxation into a feasible cyclic policy of comparable cost. For this purpose, we first reduce the prescribed run frequencies by a factor
of $2$ and round their corresponding intervals upward to powers
of $2$, obtaining a sequence of periodic due dates within a common cycle. Here, each production lot exactly matches the demand between two successive
due dates, and its setup and uninterrupted production run
form a single block to be placed on the timeline.
We schedule these blocks in reverse time,
giving priority according to each commodity's holding
cost per unit of normalized inventory, weighted by the fraction of its block spent in production.
To control holding costs, we examine each prefix of this priority
order and bound its excess occupied machine time, relative to
its average allocation, in terms of blocks outside that prefix.
The critical comparison utilizes the relaxation's pairwise inventory
bounds to charge these interruptions to inventories inside
the prefix, ultimately yielding Theorem~\ref{thm:main}.

\paragraph{Answering activity queries.} In
Section~\ref{sec:execute}, we show how to recover individual
activities from our compact description, with a running time
that depends logarithmically on the period length, measured in integer
time units.
Reading the policy construction backward turns due dates into release
times, which can be interpreted as repeatedly selecting the highest-priority
available job and processing it without interruption.
For each prefix of the priority order, we consider an auxiliary
machine serving only that prefix and compute the number of integer time
units during which this machine would be idle.
These counts depend only on the released amounts of work,
enabling us to avoid identifying the individual
jobs processed by the auxiliary machine.
The divisibility of the rounded power-of-$2$ periods reduces each
count to an explicit expression evaluated at only $O(n)$
candidate times.
We then scan the commodities from lowest to highest priority,
using these counts and binary searches to identify the latest
activity start at or before the queried time.

\section{The Convex Relaxation} \label{sec:lower}

In what follows, we describe our approach to formulating a convex relaxation that provides a computable lower bound on $\OPTdyn$ and guides the construction of a provably good production policy. For this purpose, in
Section~\ref{sec:cycle-account}, we show that restricting attention
to cyclic policies does not change the infimum long-run average cost and introduce auxiliary quantities that characterize such policies.
In Section~\ref{sec:relaxation}, we derive three families of constraints that must be satisfied by these quantities, due to commodities sharing a single machine and indirectly affecting each other's inventory levels. In Section~\ref{sec:inventory-lp}, we collect these constraints to formulate a convex relaxation that can be approximated to any prescribed accuracy.

\subsection{Cyclic approximation and single-period accounting} \label{sec:cycle-account}

\paragraph{Equality of dynamic and cyclic optima.} Rather than presenting upcoming developments in terms of arbitrarily structured dynamic policies, it would be convenient to focus on cyclic policies, whose specifics were discussed in Section~\ref{sec:cyclic-queries}. Moving forward, let $\Pi_{\mathrm{cyc}}$ be the class of feasible cyclic policies, with finitely many runs per period and with their initial stocks already chosen, and let $\OPTcyc = \inf_{ \pi \in \Pi_{\mathrm{cyc}} } \mathcal C(\pi)$. Lemma~\ref{lem:cyclic-density} below, whose proof is provided in Appendix~\ref{sec:cyclic-density-proof}, allows us to work with this class by establishing equivalence of the dynamic and cyclic optima. At a high level, our proof starts from any
feasible dynamic policy $\pi$ of finite cost and identifies a long
segment on which each commodity ends with nearly the same stock as
it began with. By slightly extending a single run of each commodity and
adding idle time if needed, we balance production and demand over
the resulting period without adding any setups. A small increase in
initial stocks then enables us to repeat this period indefinitely at
average cost at most $\mathcal C(\pi)+\eps$, for any prescribed $\eps>0$.

\begin{lemma} \label{lem:cyclic-density}
$\OPTcyc=\OPTdyn$.
\end{lemma} 

\paragraph{Single-period quantities.} For the remainder of this section, let us fix a feasible cyclic policy $\pi$ of finite cost, with period length $P$. We focus on a single cycle of the form ${\cal S} = [t_0, t_0+P)$, where $t_0$ is the beginning of an arbitrarily chosen setup. For every commodity $i \in [n]$, let $m_i^\pi \geq 1$ be its number of runs across ${\cal S}$ and let $\ell_{ir}^\pi$ be the production time of its $r$th run, excluding setup. As mentioned in Section~\ref{sec:cyclic-queries}, production must balance demand over ${\cal S}$, since a deficit would eventually exhaust our stock, whereas a surplus would result in infinite average holding costs, and thus, $\sum_{r \in [m_i^\pi]} \ell_{ir}^\pi = \rho_i P$. We also define the setup frequency of this commodity as $f_i^\pi = \frac{m_i^\pi}{P}$ and its mean production duration as $\mu_i^\pi = \frac{1}{m_i^\pi} \cdot \sum_{r \in [m_i^\pi]} \ell_{ir}^\pi = \frac{\rho_i}{f_i^\pi}$. 

For future accounting purposes, let us introduce $q_i^\pi(t) = \frac{I_i^\pi(t)}{p_i}$ as the inventory of commodity $i$ at time $t$, expressed
as the amount of machine time required to produce the quantity currently
in stock. To better understand how this function evolves, we can rewrite the inventory equation~\eqref{eqn:inventory_time} to obtain $q_i^\pi(t) = \frac{I_i^\pi(0)}{p_i}
+ |\mathcal{P}_i\cap[0,t)|
- \rho_i t$ for every $t\ge0$. Here,
each time unit spent producing $i$ adds one unit of normalized stock,
while demand removes $\rho_i$ inventory units per unit of time. Hence, $\frac{ \dd q_i^\pi(t) }{ \dd t } =\ind\{t\in\mathcal{P}_i\}-\rho_i$ for almost
every $t\ge0$.

\paragraph{Cost representation.} Given these quantities, we express the long-run average cost $\mathcal C(\pi)$ of this policy in a more workable format, by defining $Q_i^\pi = \frac{ 1 }{ P } \cdot \int_{\cal S} q^\pi_i(t) \dd t$ as the mean normalized inventory level of commodity $i$ over a single period, and $c_i=h_i p_i$ as its holding cost rate per normalized inventory unit. In turn, we conclude that $\mathcal C(\pi)$ is a linear function of $\{f_i^\pi \}_{i \in [n]}$ and $\{ Q_i^\pi \}_{i \in [n]}$, by observing that
\begin{align}
 \mathcal C(\pi)
 ~~&=~~ \frac{1}{P} \cdot \left(
      \sum_{i\in[n]}K_i m_i^\pi
      + \int_{\mathcal S}\sum_{i\in[n]}h_i I_i^\pi(t) \dd t
    \right) \nonumber        \\
 &=~~ \sum_{i\in[n]} \left( K_i f_i^\pi
    + \frac{c_i}{P} \cdot \int_{\mathcal S}q_i^\pi(t) \dd t \right) \nonumber \\
 &=~~ \sum_{i\in[n]}\left(K_i f_i^\pi+c_i Q_i^\pi\right), \label{eqn:represent_cost_var}
\end{align}
where the second equality holds since $f_i^\pi = \frac{m_i^\pi}{P}$, $c_i=h_i p_i$, and $q_i^\pi(t) = \frac{I_i^\pi(t)}{p_i}$. 

\subsection{Valid inequalities for cyclic policies} \label{sec:relaxation}

We now ask what the setup frequencies $\{ f_i^\pi \}_{i \in [n]}$ of a cyclic policy $\pi$ must satisfy within themselves, and what these statistics tell us
about the mean inventories $\{ Q_i^\pi \}_{i \in [n]}$ that must be carried. In the upcoming discussion, inequalities~1 and~2 account for the time available for setups and the stock that builds
up during a commodity's own production runs. The remaining
inequalities account for the stock that each commodity must carry
during setups and production runs of other commodities.

\paragraph{Inequality 1: Machine time.}
Over the segment $\mathcal S$, each commodity $i \in [n]$ requires $\rho_iP$ units
of production time and $m_i^\pi s_i$ units of setup time. Since setups and
production runs do not overlap, we have $\sum_{i \in [n]} (m_i^\pi s_i + \rho_iP) \leq P$, and therefore
\begin{equation} \label{eqn:UB_sum_sf}
\sum_{i\in[n]} s_i f_i^\pi
~~=~~ \frac{1}{P} \cdot \sum_{i\in[n]} m_i^\pi s_i ~~\leq~~ 1 - 
\sum_{i\in[n]} \rho_i ~~=~~ 1-R ~~=~~ \kappa.  
\end{equation}

\paragraph{Inequality 2: Individual inventories.} The following lemma, whose proof is provided in Appendix~\ref{app:proof_lem_individual}, relates the mean normalized stock $Q_i^\pi$ of each commodity $i \in [n]$ to its setup frequency $f_i^\pi$. Within this proof, $a_i^\pi(t)$ will indicate that commodity $i$ is produced at
time $t$. The crux of our argument resides in showing that the mean stock of $i$
during its production runs equals its mean stock over the
entire period. Focusing on the former term, during a production run of length $\ell$, stock rises from a
nonnegative level at rate $1-\rho_i$, implying that the area under the inventory
trajectory is at least $\frac{ (1-\rho_i)\ell^2 }{ 2 }$. For a fixed number of runs
and fixed total production time, the sum of these areas can be lower-bounded by the value obtained when all runs have length
$\mu_i^\pi= \frac{\rho_i}{f_i^\pi}$ and start with zero stock.

\begin{lemma}\label{lem:individual}
$Q_i^\pi \ge \frac{(1-\rho_i) \rho_i}{2f_i^\pi}$, for every commodity $i \in [n]$.
\end{lemma}

\paragraph{Additional inequalities: Inventory carried during external blocks.} For convenience, a setup for commodity $i$, together with the production run
immediately following it, will be called an {$i$-block}, with 
$a_i^{+\pi}(t)$ indicating that an $i$-block is in progress at time $t$. For each
ordered pair of commodities $i \neq j\in[n]$, let $Y_{ij}^\pi=\frac{1}{P} \cdot \int_{\mathcal S}
a_i^{+\pi}(t)q_j^\pi(t)\dd t$ be the contribution to the mean inventory $Q_j^\pi = \frac{ 1 }{ P } \cdot \int_{\cal S} q^\pi_j(t) \dd t$ of $j$ from the time spent setting up for and producing $i$. 

\begin{figure}[htbp!]
\centering
\resizebox{0.88\linewidth}{!}{%
\begin{tikzpicture}[els/base,x=1.8cm,y=1.1cm]


  \path[
    els/block,
    draw=elsOrange,
    pattern=north east lines,
    pattern color=elsOrange
  ] (0,3.75) rectangle (1.2,4.12);

  \path[els/block,draw=elsOrange,fill=elsOrange]
    (1.2,3.75) rectangle (5,4.12);

  \node[fill=white,inner sep=1.2pt]
    at (.6,3.935) {Setup};
  \node[text=white]
    at (3.1,3.935) {Production};

  \draw[els/axis] (0,4.37)--(5,4.37);
  \draw[els/axis] (0,4.27)--(0,4.47);
  \draw[els/axis] (5,4.27)--(5,4.47);

  \node[anchor=south] at (2.5,4.45)
    {$x=s_i+\ell_{ir}^{\pi}$};


  \fill[elsBlue!16]
    (0,0)--(0,2.5)--(5,0)--cycle;

  \foreach \x in {0,1.2,5}
    \draw[els/guide,densely dotted]
      (\x,0)--(\x,3.45);

  \draw[els/axis]
    (-.4,3.45)--(-.4,0)--(5.35,0);

  \draw[elsBlue,dashed,line width=.75pt]
    (0,2.5)--(5,0);

  \draw[els/curve,elsBlue]
    (0,3.15)--(5,.65);
  \fill[elsBlue] (5,.65) circle (1.5pt);

  \node[anchor=east] at (-.7,1.6) {$q_j^\pi(t)$};
  \node[anchor=east] at (-.45,0) {$0$};
  \node[anchor=west] at (5.40,0) {$t$};

  \foreach \x/\lab in {0/v,1.2/{v+s_i},5/{v+x}} {
    \draw[els/axis] (\x,0)--(\x,-.06);
    \node[anchor=north] at (\x,-.10) {$\lab$};
  }

  \node at (2.9,2.65) {$-\rho_j$};
  \draw[els/leader,elsBlue]
    (2.85,2.45)--(2.7,1.80);

  \node[text=elsBlue] at (1.65,.72)
    {$\displaystyle\frac{\rho_jx^2}{2}$};

\end{tikzpicture}%
}
\caption{Inventory of commodity $j\ne i$ during an $i$-block
of length $x$, including setup $s_i$ and production time
$\ell_{ir}^{\pi}$. The solid line is the actual normalized
inventory, which decreases at rate $\rho_j$. Nonnegativity
implies the dashed lower bound, whose induced shaded area is $\frac{ \rho_jx^2 }{ 2 }$.}
\label{fig:external-block-inventory}
\end{figure}
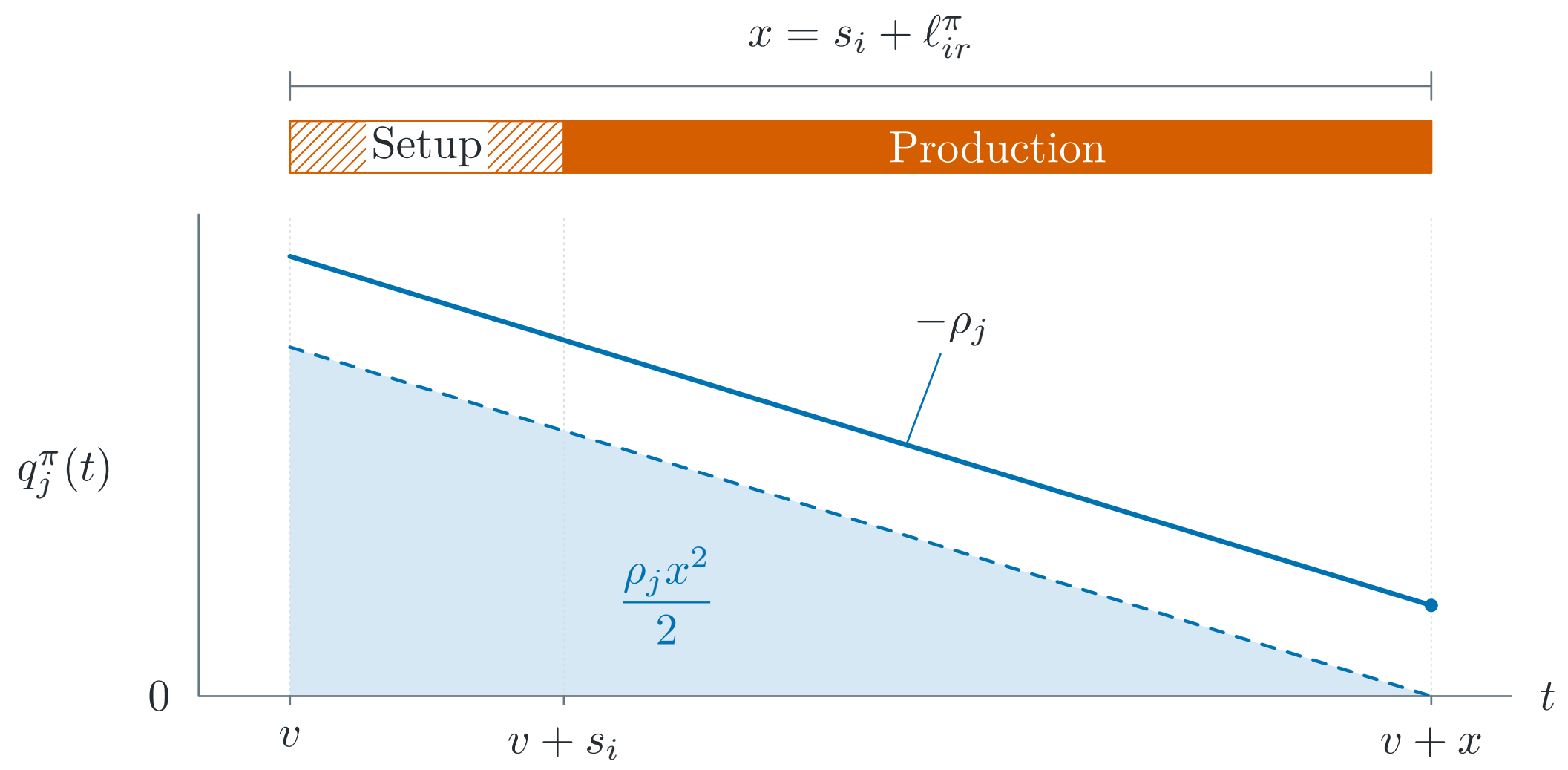

In Appendix~\ref{app:proof_lem_block-inventory}, we prove three inequalities that relate $\{  Y_{ij}^\pi \}_{i \neq j \in [n]}$ to run frequencies and mean inventories, as formally stated in Lemma~\ref{lem:block-inventory}. First, blocks of commodities other than $j$ do not overlap and can
account only for stock carried while $j$ is not being produced;
those times contribute at most 
$(1-\rho_j)Q_j^\pi$ to its mean inventory, corresponding to item~1 below. Second, by arguing that the product
$q_i^\pi q_j^\pi$ has zero net change over a period, the contributions
from carrying each commodity's stock during production of the other
sum to $\rho_iQ_j^\pi+\rho_jQ_i^\pi$, and including setups can only
increase this sum, as formalized in item~2. Third, during an $i$-block of length $x$, commodity
$j$ must meet demand entirely from stock. Hence, its normalized inventory
must cover the remaining demand at every point in this block, forcing
an area of at least $\frac{ \rho_j x^2 }{ 2 }$, and collecting these quantities leads to item~3. Figure~\ref{fig:external-block-inventory} illustrates this requirement:
The shaded triangle is the minimum inventory area compatible with
nonnegative stock throughout the block.

\begin{lemma} \label{lem:block-inventory}
The following inequalities hold:
\begin{enumerate}
\item $\sum_{i\ne j}Y_{ij}^\pi \le(1-\rho_j)Q_j^\pi$,       for every $j\in[n]$.

\item $Y_{ij}^\pi +Y_{ji}^\pi\ge\rho_iQ_j^\pi+\rho_jQ_i^\pi$,
      for every $i \neq j\in[n]$.
      
\item $Y_{ij}^\pi\ge\frac{\rho_j}{2}
      (\frac{\rho_i^2}{f_i^\pi}+2\rho_i s_i+f_i^\pi s_i^2)$,
      for every $i \neq j\in[n]$.
\end{enumerate}
\end{lemma}

\subsection{The convex program} \label{sec:inventory-lp}

We now treat the frequencies $\{ f_i \}_{i \in [n]}$ and the block contributions
$\{ Y_{ij} \}_{ i\ne j \in [n]}$ as decision variables within our convex relaxation. To simplify upcoming notation, for every commodity $j \in [n]$, we define the auxiliary function
\begin{equation} \label{eq:inventory-charge}
 Q_j(f,Y) ~~=~~
 \max\left\{
\frac{\rho_j(1-\rho_j)}{2f_j},
\frac{1}{1-\rho_j} \cdot \sum_{i\ne j}Y_{ij}
\right\}.
\end{equation}
These two terms are precisely the lower bounds
on the mean normalized stock of $j$ supplied by Lemma~\ref{lem:individual}
and item~1 of Lemma~\ref{lem:block-inventory}, respectively. By mimicking representation~\eqref{eqn:represent_cost_var} of the long-run average cost, our convex program minimizes the resulting setup and
holding costs, subject to the family of valid inequalities described in Section~\ref{sec:relaxation}. This formulation can be compactly written as follows: 
\begin{align}
 \min \quad
&\sum_{i\in[n]}(K_i f_i+c_iQ_i(f,Y)) \tag{CVX}
 \label{eq:inventory-convex}\\
 \text{s.t.}\quad
 &\sum_{i\in[n]}s_i f_i\le\kappa
 \label{eq:convex-machine}\\
 &Y_{ij}+Y_{ji}
     \ge\rho_iQ_j(f,Y)+\rho_jQ_i(f,Y)
 && \forall \, i \neq j\in[n]
 \label{eq:convex-pair}\\
 &Y_{ij}\ge\frac{\rho_j}{2}
       \left(\frac{\rho_i^2}{f_i}+2\rho_i s_i+f_i s_i^2\right)
 &&\forall \, i \neq j\in[n] \label{eq:convex-block} \\
 & f>0,Y\ge0 \nonumber
\end{align}

\paragraph{Comparison with the dynamic optimum.} 
We proceed by proving that the frequencies $\{ f_i^\pi \}_{i \in [n]}$ and block contributions
$\{ Y_{ij}^\pi \}_{ i\ne j \in [n]}$ of every cyclic policy $\pi$ jointly form a
feasible point $(f^\pi,Y^\pi)$ to formulation~\eqref{eq:inventory-convex}. This claim is embedded within the next result, whose formal proof is presented in 
Appendix~\ref{app:proof_lem_convex-lower-bound}, showing that $\OPT\eqref{eq:inventory-convex}$ indeed provides a lower bound on the dynamic optimum. 

\begin{lemma} \label{lem:convex-lower-bound}
$\OPT\eqref{eq:inventory-convex} \le \OPTdyn$.
\end{lemma}

\paragraph{Approximating \eqref{eq:inventory-convex}.} Toward obtaining a near-optimal solution to our relaxation, we first explain why \eqref{eq:inventory-convex} is a convex program, noting that it consists of  $O( n^2 )$ decision variables and $O( n^2 )$ constraints:
\begin{itemize}
    \item {\em Objective function}: The first expression in definition~\eqref{eq:inventory-charge} is convex in $f_j>0$ and the second is linear in $Y$. Their maximum $Q_j(f,Y)$ is therefore convex in $(f,Y)$. Since  $c_i>0$ for every commodity $i \in [n]$, the objective function $\sum_{i\in[n]}(K_i f_i+c_iQ_i(f,Y))$ is convex as well.

    \item {\em Constraints}: The machine-time constraint~\eqref{eq:convex-machine} is linear. For the pair constraint~\eqref{eq:convex-pair}, the right-hand side is a positive weighted sum of $Q_j(f,Y)$ and $Q_i(f,Y)$, which are convex, as explained above.  For the block constraint~\eqref{eq:convex-block}, its right-hand side is a positive multiple of $\frac{ 1 }{ f_i }$ plus an affine function of $f_i$. In both cases, the left-hand side is linear, meaning that each inequality defines a convex set. 
\end{itemize}

Given these properties, for every prescribed accuracy $\eps > 0$, we can compute a rational feasible solution to \eqref{eq:inventory-convex} whose objective value is within a factor of $1+\eps$ of optimal,
in deterministic polynomial time in the binary encoding lengths of the instance considered and $\frac{ 1 }{ \eps }$. This result follows from standard convex optimization methods \citep[Thm.~6.3.7]{Lovasz2003SDP}, since our program admits a polynomial-time separation oracle, and since an
optimum-preserving cap on its objective creates a bounded epigraph with polynomially encoded inner and outer ball bounds. Moreover, a computable positive rational lower bound on the optimum can be utilized to convert additive accuracy into the stated multiplicative guarantee. For completeness, Appendix~\ref{app:convex-solvability}
provides a short verification of these conditions. 
\section{Constructing and Analyzing a Cyclic Policy} \label{sec:construction}

In what follows, we explain how an approximate solution of the convex
relaxation can guide the construction of a feasible cyclic policy whose long-run average cost
is within a constant factor of $\OPT\eqref{eq:inventory-convex}$. While our relaxation specifies how frequently each commodity should be produced, it
does not provide direct information on when its runs should take place. In
Section~\ref{sec:calendar-data}, we utilize these frequencies to choose
lot sizes, which are coupled with due dates across a common period.
In Section~\ref{sec:calendar-rule}, we schedule the resulting lots in reverse time, based on a priority rule that chooses among those
that can be completed at the current time.
In Section~\ref{sec:calendar-inventories}, we determine initial stocks
that allow this schedule to be repeated while meeting demand at every
point and bound the setup cost incurred. The crux of our analysis is revealed in Section~\ref{sec:occupation}, which relies on the specific priority rule in play to bound the
policy's holding cost in terms of $\OPT\eqref{eq:inventory-convex}$.

\subsection{The auxiliary scheduling problem} \label{sec:calendar-data}

\paragraph{Approximate solution to~\eqref{eq:inventory-convex}.} Following the discussion in Section~\ref{sec:inventory-lp}, for any fixed accuracy level $\eps > 0$, we can
compute in deterministic polynomial time a rational feasible point
$({f}, {Y})$ whose objective value in~\eqref{eq:inventory-convex} is within a factor of $1 + \eps$ of optimal. Letting ${Q}_i=Q_i( {f}, {Y})$
be the inventory charges in definition~\eqref{eq:inventory-charge}, we introduce the shorthand notation ${S} = \sum_{i \in [n]} K_i {f}_i$ and ${H} = \sum_{i \in [n]} c_i {Q}_i$ for the setup and holding terms of the relaxation objective, meaning that ${S} + {H} \le (1+\eps) \cdot \OPT\eqref{eq:inventory-convex}$. Moving forward, it will be convenient to pick $\eps = \frac{ 1 }{ 120 }$; this accuracy will suffice for our final approximation guarantee. This way, by defining $L  = \frac{ 120 }{ 121 } \cdot (S+H)$, it follows that $L \leq \OPT\eqref{eq:inventory-convex} \leq \OPTdyn$. Finally, in analogy to Section~\ref{sec:relaxation}, we refer to $\mu_i= \frac{ \rho_i }{ {f}_i }$ as the mean production
duration of each commodity $i \in [n]$, noting that at frequency ${f}_i$, runs of duration $\mu_i$ would produce at the average rate
$p_i {f}_i \mu_i=d_i$, exactly balancing demand.

\paragraph{Choosing lot sizes and due-date intervals.} Given $({f}, {Y})$, for every commodity $i \in [n]$, we define its individual period length $T_i$ by rounding $\frac{2}{{f}_i}$ upward to the nearest power of $2$, meaning that $T_i = 2^{\lceil\log_2(2/{f}_i)\rceil}$. This choice ensures that the periods $\{ T_i\}_{i \in [n]}$ have a common arithmetic structure, noting that each resulting period has polynomial encoding length, with
\begin{equation}\label{eq:dilation}
\frac{ 2 }{ f_i } ~~\le~~ T_i ~~<~~ \frac{ 4 }{ f_i }.
\end{equation}
Next, with one lot due every $T_i$ time units, each lot must produce $d_iT_i$ units of commodity $i$. Its production duration is therefore set to $\ell_i = \frac{d_iT_i}{p_i} = {f}_iT_i \mu_i$, meaning that the combined length of a setup-production block becomes $e_i = s_i + \ell_i$.

For purposes of analysis, it is instructive to associate each commodity with its utilization statistic, $u_i$. This parameter represents the fraction of machine time occupied by commodity $i$ when one block of duration $e_i$ is scheduled every $T_i$ time units on average, i.e., $u_i=\frac{e_i}{T_i}=\rho_i+\frac{s_i}{T_i}$. Letting $U=\sum_{i\in[n]}u_i$ designate total utilization, a key property that will be exploited by our construction is that 
\begin{equation} \label{eq:calendar-load}
 U ~~=~~ \sum_{i\in[n]} \left( \rho_i+\frac{s_i}{T_i} \right)
 ~~\le~~ R+\frac{1}{2} \cdot \sum_{i\in[n]} s_i{f}_i
 ~~\le~~ 1 - \kappa +\frac{\kappa}{2}
 ~~<~~1.
\end{equation}
Here, the first inequality holds since $T_i \geq \frac{ 2 }{ {f}_i }$ by relation~\eqref{eq:dilation}. The second inequality is obtained by recalling that $1-R = \kappa$ and $\sum_{i\in[n]}s_i {f}_i\le\kappa$, by constraint~\eqref{eq:convex-machine}, as $({f},{Y})$ is feasible in~\eqref{eq:inventory-convex}. 

\paragraph{The scheduling problem.} Having ensured that $\{ T_i \}_{ i \in [n]}$ are all powers of $2$, our joint period length will be $P=\max_{i\in[n]}T_i$. Over this period, we wish to schedule $N_i=\frac{P}{T_i}$ lots of each commodity $i \in [n]$, noting that $N_i$ is a positive
integer. The period $P$ and the counts $\{ N_i \}_{ i \in [n]}$ can be computed and
represented in polynomial time. The auxiliary scheduling problem we consider consists of $N_i$ lots of each commodity
$i$, indexed by $r \in [N_i]$. Lot $(i,r)$ requires a
block of length $e_i$ and its due date is $d_{ir}=rT_i$. We aim to construct a
schedule of these blocks within $[0,P)$, avoiding overlaps,
such that each completion time $C_{ir}$ does not exceed the corresponding due date
$d_{ir}$. A block completed at time  $C_{ir}$ consists of
setup on $[C_{ir}-e_i,C_{ir}-\ell_i)$ followed by production
on $[C_{ir}-\ell_i,C_{ir})$. Lots may be completed before their due dates, and
idle intervals are allowed. 

\subsection{Construction of a priority-based schedule} \label{sec:calendar-rule}

\paragraph{Existence of a feasible schedule.}
This auxiliary setting is a classical single-machine scheduling
problem, where all jobs are available at time $0$, with specified processing times and due dates. In this context, 
Jackson's earliest-due-date (EDD) rule schedules jobs in nondecreasing
order of their due dates and produces a feasible schedule
whenever one exists 
\citep{Jackson1955,GrahamEtAl1979}. For our specific instance, feasibility directly follows  from the load
bound~\eqref{eq:calendar-load}. Indeed, for every $0<t\le P$,
the total length of all blocks due by time $t$ is
\begin{equation} \label{eq:deadline-work}
 \sum_{i\in[n]} e_i \cdot        \left\lfloor\frac{t}{T_i}\right\rfloor
 ~~\le~~ t \cdot \sum_{i\in[n]}\frac{e_i}{T_i}
 ~~=~~ Ut ~~<~~ t.
\end{equation}
Hence, if all lots are processed consecutively from
time $0$ in earliest-due-date order, the completion time of
each lot is at most the total length of all blocks due no
later than that lot. Applying inequality~\eqref{eq:deadline-work} at its
due date shows that it is completed on time.

\paragraph{Controlling holding costs: Feasibility is insufficient.} That said, our objective is to construct a feasible schedule whose
induced holding cost can be bounded in terms of $\OPT\eqref{eq:inventory-convex}$.
The timing of production matters for this purpose, since completing
a lot too early may require carrying its output for longer before demand consumes it. As explained below, our construction works backward, from time $P$ to time $0$. At each step, we choose a lot whose
due date permits it to be completed at the current time, with specific priorities determining the choice among eligible commodities to allow subsequent comparisons with $\OPT\eqref{eq:inventory-convex}$.

\paragraph{Priorities and the placement rule.} Our construction assigns all lots of each commodity $i \in [n]$ a uniform priority of $w_i=\frac{c_i \ell_i}{e_i}$, where $c_i$ is the holding cost coefficient in the objective function of~\eqref{eq:inventory-convex}, and $\frac{\ell_i}{e_i}$ is the proportion of an $i$-block spent in production. We initialize the clock at $t=P$, and let $m_i=N_i$ be the number
of still-unplaced lots of commodity $i$. Lots of this commodity will be placed in decreasing order of their indices, implying that the remaining
due dates are $T_i,2T_i,\ldots,m_iT_i$. Now, as long as some lot remains, we employ the following rule:
\begin{itemize}
    \item If $m_i \geq 1$ and $m_iT_i\ge t$ for at least one commodity $i$:
    \begin{itemize}
        \item Choose a commodity $i$ of highest priority $w_i$ among those satisfying these conditions, breaking ties in favor of the smaller index.

        \item Place lot $(i,m_i)$ in $[t-e_i,t)$, with its setup first and its production run second.

        \item Set $m_i\leftarrow m_i-1$ and $t\leftarrow t-e_i$.
    \end{itemize}  

    \item Otherwise, let $d=\max_{i:m_i \geq 1} \{ m_iT_i \}$. Leave the interval $[d,t)$ idle and set $t\leftarrow d$.
\end{itemize}

\paragraph{Feasibility of the backward placement rule.}
The eligibility condition proposed above ensures that each lot finishes by its due
date, assuming it is actually placed. However, we must still verify that every lot is indeed placed without running out of time, particularly when idle intervals are inserted and the clock is moved backward. The following claim, whose proof is provided in Appendix~\ref{app:proof_lem_calendar-fit}, shows that similarly to Jackson's earliest-due-date rule, we are successful in placing all lots. The benefits of our rule will materialize in the next section, when we convert the resulting schedule into a production policy.

\begin{lemma} \label{lem:calendar-fit}
The placement rule schedules all lots in nonoverlapping blocks
within $[0,P)$. Moreover, $C_{ir}\le d_{ir}$ for every lot
$(i,r)$.
\end{lemma}

\subsection{Translating the schedule to a cyclic policy}
\label{sec:calendar-inventories}

\paragraph{The resulting policy and feasibility.} Let ${\pi}$ be the cyclic policy obtained by replicating
the schedule constructed in Section~\ref{sec:calendar-rule} every $P$ time units, including its setups, production runs, and idle intervals. As the initial inventory level of each commodity $i \in [n]$, we set $I_i^{{\pi}}(0) = p_i(1-\rho_i)\ell_i$. These stocks are finite rational numbers of polynomial
encoding length.

We proceed by proving that, for each commodity, the chosen initial stock suffices to establish that demand is satisfied at any point in time, even if every production run ends exactly at its due date. Moreover,
production balances demand over each period, implying that feasibility
during the first period extends to the entire planning horizon.

\begin{lemma}\label{lem:calendar-feasible}
The policy ${\pi}$ is feasible. Moreover,
$I_i^{{\pi}}(t+P)=I_i^{{\pi}}(t)$ for every
commodity $i\in[n]$ and time $t\ge0$.
\end{lemma}
\begin{proof}
Focusing on a single commodity $i \in [n]$, lot $(i,r)$ is produced on $[C_{ir}- \ell_i,C_{ir})$, and  its contribution to cumulative
production time by time $t$ is
$\min\{\ell_i, [t-C_{ir}+\ell_i]^+\}$.
Consequently, for every $t\in[0,P]$, due to starting with an inventory level of $I_i^{{\pi}}(0) = p_i(1-\rho_i)\ell_i$, we have
\[
 q_i^{{\pi}}(t) ~~=~~ \frac{ I_i^{{\pi}}(t) }{ p_i }
 ~~=~~ (1-\rho_i)\ell_i-\rho_i t
       +\sum_{r \in [N_i]}
          \min\{\ell_i,[t-C_{ir}+\ell_i]^+\}.
\]

We claim that $q_i^{{\pi}}(t) \geq 0$ for every $t\in[0,P)$. To this end, let us write $t=mT_i+u$, where $m = \lfloor \frac{ t }{ T_i } \rfloor$ and $0\le u<T_i$. By Lemma~\ref{lem:calendar-fit}, the first $m$ lots of commodity
$i$ have already been completed, since
$C_{ir}\le d_{ir} = rT_i\le t$ for every $r\le m$. The next lot completes by $(m+1)T_i$, and its contribution to
production time by $t$ is at least
$[ u-(T_i-\ell_i) ]^+$. Since all remaining contributions are nonnegative,
\begin{align*}
 q_i^{{\pi}}(t)
 ~~&\ge~~
 (1-\rho_i)\ell_i-\rho_i(mT_i+u)
 +m\ell_i+[ u-(T_i-\ell_i) ]^+\\
 &=~~ (1-\rho_i)\ell_i-\rho_i u
    +[ u-(T_i-\ell_i) ]^+\\
 &=~~
 \begin{cases}
  \rho_i(T_i-\ell_i-u),
     & \text{if } 0\le u\le T_i-\ell_i\\[2pt]
  (1-\rho_i)(u-T_i+\ell_i),
     & \text{if }T_i-\ell_i<u
 \end{cases}
\end{align*}
where both equalities follow by substituting $\ell_i=\rho_iT_i$. Both expressions in
the last line are nonnegative, proving that stockouts do not occur during the first period, $[0,P)$.

Finally, the total production time of commodity $i$ in a single period
is $N_i\ell_i=\rho_iP$. Hence, $q_i^{{\pi}}(P) = q_i^{{\pi}}(0)+N_i\ell_i-\rho_iP
= q_i^{{\pi}}(0)$. Since our schedule repeats every $P$ time units, its inventory trajectory repeats as well. Thus, feasibility holds for all
$t\ge0$, and multiplying by $p_i$ yields the claimed periodicity, $I_i^{{\pi}}(t+P)=I_i^{{\pi}}(t)$.
\end{proof}

\paragraph{Cost analysis.} Let us recall that, by representation~\eqref{eqn:represent_cost_var}, the long-run average cost of the policy ${\pi}$ can be decomposed into its setup and holding costs, designated by ${\cal C}_{\mysetup}({\pi}) = \sum_{i\in[n]} K_i f_i^{{\pi}}$ and ${\cal C}_{\myhold}( {\pi} ) = \sum_{i\in[n]} c_i Q_i^{{\pi}}$. These two components are accounted for as follows:
\begin{itemize}
    \item {\em Setup cost}: For every commodity $i$, the policy ${\pi}$ performs exactly $N_i$ setups in each period of length $P$, leading to a setup frequency of $\frac{N_i}{P}=\frac{1}{T_i}$. Consequently, its long-run average setup cost is
    \[ {\cal C}_{\mysetup}({\pi})  ~~=~~ \frac{1}{P}\cdot\sum_{i\in[n]}K_iN_i  ~~=~~ \sum_{i\in[n]}\frac{K_i}{T_i}  ~~\le~~ \frac{1}{2}\cdot\sum_{i\in[n]}K_i{f}_i  ~~=~~ \frac{{S}}{2}, \]
    where the inequality holds since $T_i \ge \frac{ 2 }{ {f}_i }$, by relation~\eqref{eq:dilation}.

    \item {\em Holding cost}: Bounding holding costs is the central difficulty in our analysis, since meeting due dates alone does not control the amount of inventory that may accumulate when lots are produced early. Our priority rule allows us to bound how much the time devoted to higher-priority commodities exceeds their average allocation, while the relaxation's inventory constraints provide the quantities against which this excess can be charged. Section~\ref{sec:occupation} develops this connection to establish the following upper bound on the long-run average holding cost of ${\pi}$.    

    \begin{lemma} \label{lem:calendar-cost}
    ${\cal C}_{\myhold}( {\pi} ) \leq 4(3+\sqrt2)\cdot {H}$.
    \end{lemma}
\end{itemize}
Synthesizing these two bounds, we infer the approximation guarantee stated in Theorem~\ref{thm:main} by observing that the long-run average cost of $\pi$ is
\begin{align}
{\cal C}({\pi}) ~~&=~~ {\cal C}_{\mysetup}({\pi}) + {\cal C}_{\myhold}( {\pi} ) \nonumber \\
&\le~~ \frac{S}{2} + 4(3+\sqrt2) \cdot H \nonumber \\
~~&\le~~ 4(3+\sqrt2) \cdot (S+H) \nonumber \\
&\leq~~ \frac{121(3+\sqrt2)}{30} \cdot \OPT\eqref{eq:inventory-convex} \label{eqn:final_approx_1} \\
&\leq~~ 18 \cdot \OPTdyn. \label{eqn:final_approx_2}   
\end{align}
Here, inequality~\eqref{eqn:final_approx_1} holds since 
$\frac{ 120 }{ 121 } \cdot (S+H) \leq \OPT\eqref{eq:inventory-convex}$, as explained in Section~\ref{sec:calendar-data}, and inequality~\eqref{eqn:final_approx_2} follows by recalling that $\OPT\eqref{eq:inventory-convex} \le \OPTdyn$, according to Lemma~\ref{lem:convex-lower-bound}.

\subsection{Proof of Lemma~\ref{lem:calendar-cost}: Holding cost bound} \label{sec:occupation}

\paragraph{Inventory and priority prefixes.}

For purposes of analysis, we relabel commodities in weakly decreasing order
of priority, $w_1\ge\cdots\ge w_n$, preserving the tie-breaking order in
Section~\ref{sec:calendar-rule}. In addition, let
$v_k=\sum_{i\in[k]}u_i$ be the sum of utilization parameters over the first $k$ commodities in this order. Recalling that
$a_i^{+{\pi}}(s)$ indicates whether an $i$-block is in
progress at time $s$, we define $D_k(t) = \int_0^t \sum_{i\in[k]} a_i^{+{\pi}}(s)\dd s - v_k t$ for every $t\in[0,P]$. Here, the first term captures the time occupied by blocks of the first $k$ commodities
during $[0,t)$, whereas $v_k t$ is the time that would have been occupied had they received their average share of the machine
throughout this interval. Letting
$\bar{D}_k=\frac{1}{P} \cdot \int_0^P D_k(t) \dd t$ be its mean  over the entire period, we proceed by bounding this mean in terms of the inventory charges $\{ Q_j \}_{j \in [k]}$, arguing that
\begin{equation} \label{eqn:prefix}
\bar{D}_k ~~\le~~ 4(1+\sqrt2) \cdot \sum_{j \in [k]} \frac{u_j Q_j}{\rho_j}. 
\end{equation}

To this end, let us refer to commodities $1, \ldots, k$ as the prefix. Over one period, prefix blocks occupy a fraction
$v_k=\sum_{j\in[k]}u_j$ of the machine's time. The fraction spent
producing prefix commodities is $r_k=\sum_{j\in[k]}\rho_j$, and
the fraction spent on their setups is $v_k-r_k$, noting that $0<r_k\le v_k\le U<1$. With this notation, step~1 below bounds $\bar{D}_k$ in terms of blocks {outside} the prefix. Step~2 then utilizes our relaxation to show how those blocks give rise to inventory charges {inside} the prefix.

\begin{figure}[htbp!]
\centering
\resizebox{0.88\linewidth}{!}{%
\begin{tikzpicture}[els/base,x=1.15cm,y=1.8cm]


  \path[els/block,draw=elsBlue,fill=elsBlue]
    (.8,2.25) rectangle (2.8,2.58);

  \path[els/block,draw=elsBlue,fill=elsBlue]
    (2.8,2.25) rectangle (4,2.58);

  \draw[white,line width=.8pt]
    (2.8,2.25)--(2.8,2.58);

  \path[els/block,draw=elsOrange,fill=elsOrange]
    (4,2.25) rectangle (8.5,2.58);

  \node[text=white] at (1.8,2.415) {$j\le k$};
  \node[text=white] at (3.4,2.415) {$j\le k$};
  \node[text=white] at (6.25,2.415) {$i>k$};

  \draw[els/axis] (4,2.84)--(8.5,2.84);
  \draw[els/axis] (4,2.74)--(4,2.94);
  \draw[els/axis] (8.5,2.74)--(8.5,2.94);

  \node[anchor=south] at (6.25,2.91) {$e_i$};


  \fill[elsBlue!16]
    (2,0)--(4,1.2)--(4,0)--cycle;

  \fill[elsOrange!14]
    (4,0)--(4,1.2)--(7,0)--cycle;

  \foreach \x in {2,4,7,8.5}
    \draw[els/guide,densely dotted]
      (\x,-.83)--(\x,1.75);

  \draw[els/axis] (.5,-.83)--(.5,1.75);
  \draw[els/axis] (.5,0)--(9,0);

  \draw[els/curve,elsMid]
    (.8,-.72)--(2,0);

  \draw[els/curve,elsBlue]
    (2,0)--(4,1.2);

  \draw[els/curve,elsOrange]
    (4,1.2)--(7,0);

  \draw[els/curve,elsMid]
    (7,0)--(8.5,-.6);

  \foreach \x/\y in {2/0,4/1.2,7/0,8.5/-.6}
    \fill[elsInk] (\x,\y) circle (1.5pt);

  \node[anchor=east] at (.15,.60) {$D_k(t)$};
  \node[anchor=east] at (.43,0) {$0$};
  \node[anchor=west] at (9.05,0) {$t$};

  \foreach \x/\lab in {2/x,4/z,7/y,8.5/b}
    \node[anchor=north] at (\x,-.90) {$\lab$};

  \node[anchor=south west] at (4.08,1.28) {$D_k(z)$};

  \node at (1.70,1.38) {$1-v_k$};
  \draw[els/leader,elsBlue]
    (1.95,1.22)--(3,.6);

  \node at (6.20,1.32) {$-v_k$};
  \draw[els/leader,elsOrange]
    (6.05,1.15)--(5.5,.6);

\end{tikzpicture}%
}
\caption{A positive interval $(x,y)$ of the prefix excess $D_k$.
The rising portion spans prefix blocks, whereas the falling
portion lies within one outside block $[z,b)$. Its slopes are
$1-v_k$ and $-v_k$, respectively. Since $y-z\le e_i$, the shaded
area is at most $\frac{ v_ke_i^2 }{ 2(1-v_k) }$.}
\label{fig:prefix-excess-triangle}
\end{figure}
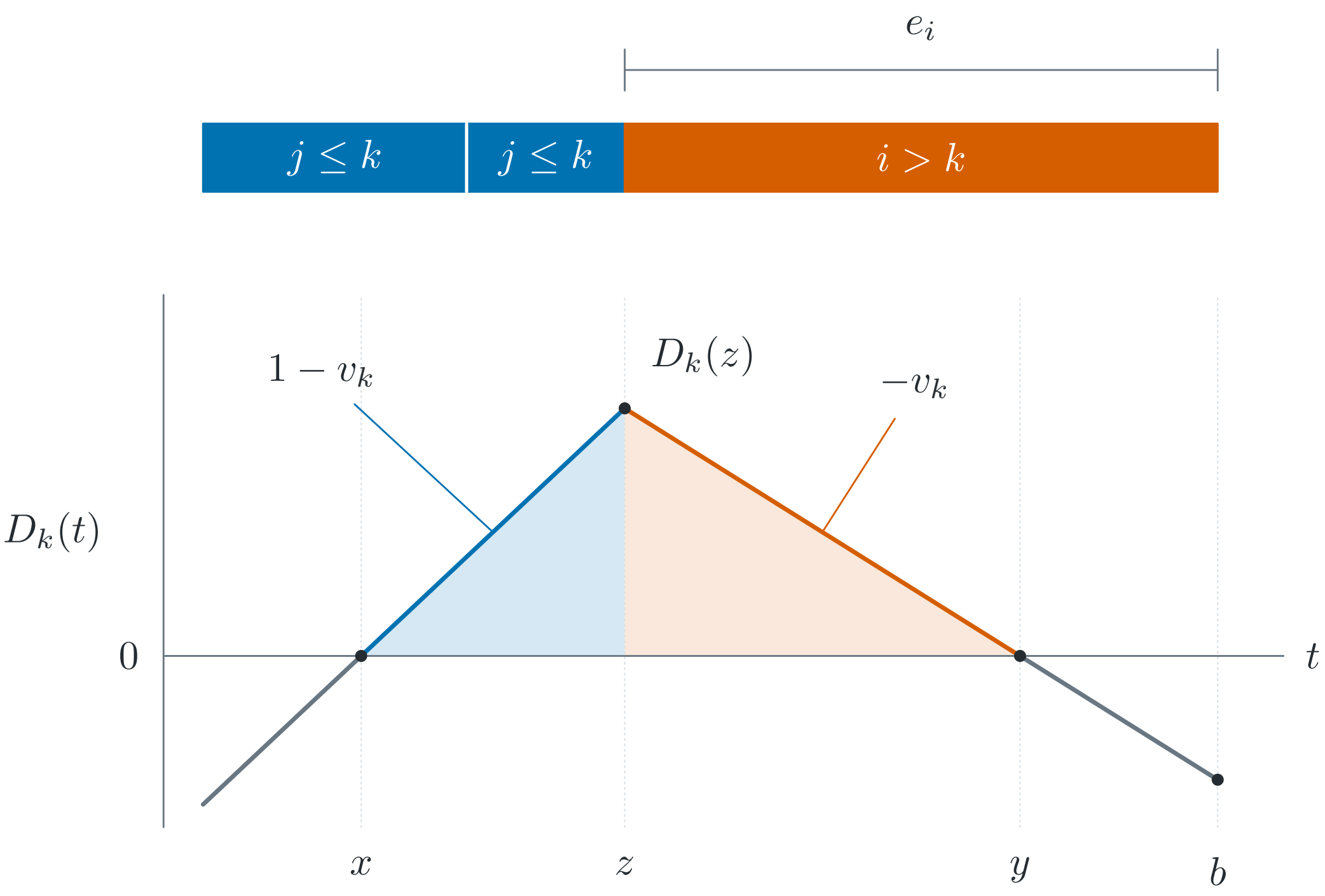

\paragraph{Step 1: Bounding the prefix's excess time.} An uninterrupted block outside the prefix can force prefix lots
to be completed earlier to meet their due dates, increasing the
machine time allocated to the prefix before that block.
To quantify these interruptions, let us  introduce
$E=\frac12 \cdot \sum_{i>k}f_i(s_i+\mu_i)^2$, which weights the squared nominal
lengths of outside blocks by their frequencies.
The following claim bounds the mean excess $\bar{D}_k$ in terms
of this quantity and the prefix utilization $v_k$; its proof is
provided in Appendix~\ref{app:proof_clm_Dk_first_UB}. The core geometric argument behind this bound is illustrated in
Figure~\ref{fig:prefix-excess-triangle}. On each maximal interval
where $D_k$ is positive, it rises while the machine serves prefix commodities and returns to zero within a single block outside the prefix.

\begin{claim}\label{clm:Dk_first_UB}
$\bar{D}_k\le \frac{ 4v_k }{ 1-v_k} \cdot E$.
\end{claim}

\paragraph{Step 2: Relation to inventory inside the prefix.}
During a block outside the prefix, every prefix commodity must
meet demand from stock, with longer interruptions requiring
more inventory to be carried for longer.
The relaxation's inventory constraints capture this requirement,
allowing the quantity $E$ from step~1 to be bounded by the inventory charges of commodities
inside the prefix.
The following claim, whose proof is
presented in Appendix~\ref{ap:proof_clm_prefix-relaxation}, formalizes this connection.

\begin{claim}\label{clm:prefix-relaxation}
$E \le \frac{ 1-r_k }{ r_k^2 - r_k v_k + v_k } \cdot \sum_{j\in[k]}\frac{u_j Q_j}{\rho_j}$.
\end{claim}

\paragraph{Deriving inequality~\eqref{eqn:prefix}.} By combining Claims~\ref{clm:Dk_first_UB} and~\ref{clm:prefix-relaxation}, we conclude that
\begin{align*}
\bar{D}_k ~~&\le~~ \frac{4v_k}{1-v_k} \cdot E\\
&\le~~ \frac{4v_k}{1-v_k} \cdot \frac{ 1-r_k }{ r_k^2 - r_k v_k + v_k } \cdot \sum_{j\in[k]}\frac{u_j Q_j}{\rho_j} \\
&\le~~ 4(1+\sqrt{2}) \cdot \sum_{j\in[k]}\frac{u_j Q_j}{\rho_j},     
\end{align*}
where the last inequality follows from the next auxiliary claim, whose proof is provided in Appendix~\ref{app:proof_clm_prefix-load}.

\begin{claim}\label{clm:prefix-load}
$\frac{v_k}{1-v_k} \cdot \frac{ 1-r_k }{ r_k^2 - r_k v_k + v_k } \le 1+\sqrt{2}$. 
\end{claim}

\paragraph{Relating ${\cal C}_{\myhold}({\pi})$ to $\{ \bar{D}_k \}_{k \in [n]}$.} The following lemma bounds the long-run average holding cost of $\pi$ through a term determined
by our initial inventories and a weighted sum of the mean prefix
excesses, to which we can apply previously established bounds.

\begin{lemma} \label{lem:chold_vs_Dk}
${\cal C}_{\myhold}({\pi}) \le \sum_{i\in[n]}c_i(1-\rho_i)\ell_i        +\sum_{k \in [n]}(w_k-w_{k+1})\bar{D}_k$.    
\end{lemma}
\begin{proof}
We first observe that for every commodity $i\in[n]$, its  inventory at any time $t\in[0,P]$ can be expressed as
\begin{align}
 I_i^{{\pi}}(t)
 ~~&=~~ I_i^{{\pi}}(0)
      +p_i \cdot \int_0^t a_i^{{\pi}}(s) \dd s-d_it \nonumber \\
&\leq~~ p_i(1-\rho_i)\ell_i + \frac{ p_i\ell_i }{ e_i } \cdot \int_0^t a_i^{+{\pi}}(s) \dd s-d_it \nonumber \\
 &=~~ p_i(1-\rho_i)\ell_i
+\frac{p_i\ell_i}{e_i} \cdot \left(\int_0^t a_i^{+{\pi}}(s)\dd s-u_it\right). \label{eqn:lem_chold_vs_Dk_0}
\end{align}
The inequality above holds since $I_i^{{\pi}}(0) = p_i(1-\rho_i)\ell_i$, as determined in Section~\ref{sec:calendar-inventories}, and since the production time of commodity $i$ across $[0,t)$ forms a fraction of at most $\frac{ \ell_i }{ e_i }$ out of the time occupied by $i$-blocks. The last equality follows by recalling that $u_i=\frac{e_i}{T_i}$ and $\ell_i = \frac{d_iT_i}{p_i}$.

With this observation, the average holding cost of ${\pi}$ can be upper-bounded by noting that
\begin{align}
{\cal C}_{\myhold}({\pi}) ~~&=~~ \frac{ 1 }{ P } \cdot \int_0^P \sum_{i \in [n]} h_i I_i^{{\pi}}(t) \dd t \nonumber \\
&\leq~~ \frac{ 1 }{ P } \cdot \int_0^P \sum_{i \in [n]} h_i \cdot \left( p_i(1-\rho_i)\ell_i       +\frac{p_i\ell_i}{e_i} \cdot \left(\int_0^t a_i^{+{\pi}}(s)\dd s-u_it\right) \right) \dd t \label{eqn:lem_chold_vs_Dk_1} \\
&=~~ \sum_{i\in[n]}c_i(1-\rho_i)\ell_i
+\frac{1}{P} \cdot \int_0^P \sum_{i\in[n]}w_i \cdot 
\left(\int_0^t a_i^{+{\pi}}(s)\dd s-u_i t\right)\dd t \label{eqn:lem_chold_vs_Dk_2} \\
&=~~ \sum_{i\in[n]}c_i(1-\rho_i)\ell_i
+\frac{1}{P} \cdot \int_0^P   \sum_{k\in[n]}(w_k-w_{k+1})D_k(t)\dd t \label{eqn:lem_chold_vs_Dk_3} \\
&=~~ \sum_{i\in[n]}c_i(1-\rho_i)\ell_i    +\sum_{k\in[n]}(w_k-w_{k+1})\bar{D}_k . \label{eqn:lem_chold_vs_Dk_4}
\end{align}
Here, inequality~\eqref{eqn:lem_chold_vs_Dk_1} is obtained by plugging in~\eqref{eqn:lem_chold_vs_Dk_0}. Equality~\eqref{eqn:lem_chold_vs_Dk_2} follows by substituting $c_i=h_ip_i$
and $w_i= \frac{ c_i\ell_i }{ e_i }$. Equality~\eqref{eqn:lem_chold_vs_Dk_3} holds since the term $\int_0^t a_i^{+{\pi}}(s)\dd s-u_i t$ appears in $D_k(t)$
precisely when $k\ge i$, and its coefficient is
therefore $\sum_{k=i}^n(w_k-w_{k+1})=w_i$, with the convention that $w_{n+1}=0$. Finally, equality~\eqref{eqn:lem_chold_vs_Dk_4} follows from the definition of $\bar{D}_k$.
\end{proof}

\paragraph{Concluding the proof of Lemma~\ref{lem:calendar-cost}.} Given the upper bound supplied by Lemma~\ref{lem:chold_vs_Dk}, we are now ready to relate the long-run average holding cost of our policy to the convex relaxation~\eqref{eq:inventory-convex}. For this purpose, we observe that
\begin{align}
{\cal C}_{\myhold}({\pi}) ~~&\le~~ \sum_{i\in[n]}c_i(1-\rho_i)\ell_i        +\sum_{k \in [n]}(w_k-w_{k+1})\bar{D}_k \nonumber \\
&\leq~~ 4 \cdot \sum_{i\in[n]}c_i(1-\rho_i) \mu_i        + 4(1+\sqrt2) \cdot\sum_{k \in [n]}(w_k-w_{k+1}) \cdot  \sum_{j \in [k]} \frac{u_j Q_j}{\rho_j} \label{eqn:lem_calendar-cost_last_1} \\
&\leq~~ 8 \cdot \sum_{i\in[n]}c_i Q_i + 4(1+\sqrt2) \cdot  \sum_{j \in [n]} \frac{w_j u_j Q_j}{\rho_j} \label{eqn:lem_calendar-cost_last_2} \\
&=~~ 4(3+\sqrt2) \cdot \sum_{i\in[n]}c_i Q_i \label{eqn:lem_calendar-cost_last_3} \\
&=~~ 4(3+\sqrt2) \cdot H. \label{eqn:lem_calendar-cost_last_4}
\end{align}
Here, inequality~\eqref{eqn:lem_calendar-cost_last_1} holds since $\ell_i = {f}_iT_i \mu_i$ and $T_i < \frac{ 4 }{ f_i }$, according to relation~\eqref{eq:dilation}, and since $\bar{D}_k \le 4(1+\sqrt2) \cdot \sum_{j \in [k]} \frac{u_j Q_j}{\rho_j}$ by inequality~\eqref{eqn:prefix}. Inequality~\eqref{eqn:lem_calendar-cost_last_2} follows by noticing that definition~\eqref{eq:inventory-charge} implies $Q_i \geq \frac{\rho_i(1-\rho_i)}{2f_i} = \frac{\mu_i(1-\rho_i)}{2}$. Finally, equalities~\eqref{eqn:lem_calendar-cost_last_3} and~\eqref{eqn:lem_calendar-cost_last_4} are respectively obtained by recalling that $w_j=\frac{c_j \ell_j}{e_j}$ and ${H} = \sum_{i \in [n]} c_i {Q}_i$.
\section{Answering Activity Queries}\label{sec:execute}

Even though our cyclic policy  may contain exponentially many runs per period, we show that its activity queries at any rational time $t\ge 0$ can be efficiently answered, thereby establishing Theorem~\ref{thm:activity-queries}. Such queries are expected to return the activity in progress at time $t$ (setup, production, or
idling), its start and end times, and the
commodity involved, when applicable. For this purpose, in Section~\ref{sec:query-reverse}, we read the schedule backward, so that lot due dates become job release times, and identify the reverse tick corresponding to $t$.
In  Section~\ref{sec:query-free}, we explain how to count the number of ticks left idle before any specified tick by an auxiliary machine serving only a given subset of commodities. Finally, in
Section~\ref{sec:query-starts}, we prove that these counts can be utilized to scan the
commodities from lowest to highest priority, locate the block or
idle interval containing the queried tick, and recover the endpoints of this activity.

\subsection{Reading the schedule backward} \label{sec:query-reverse}

\paragraph{From due dates to release times.} Our placement rule constructs one period backward, from time $P$ to time $0$, meaning that when the reverse clock reads $x$, it corresponds to time $P-x$ in the original schedule. A lot with due date $d$ can be completed
at that time precisely when $P-x\le d$, or equivalently, $x\ge P-d$. Thus, in reverse time, this lot has release time $P-d$. Since the due dates of commodity $i$ are $T_i,2T_i,\ldots,P$, its jobs are released at times $0,T_i,\ldots,P-T_i$. From this perspective, whenever the machine becomes available, it selects the highest-priority waiting job and processes it to completion;  ties between commodities are broken as in
Section~\ref{sec:calendar-rule}. Within each commodity, jobs are taken in release order, since the earliest remaining release corresponds to
the latest remaining due date in the original schedule. When no job is waiting, the machine idles until the next release.

\paragraph{Translating a physical time to a reverse tick.} To express the basic time units behind this construction, let $\Delta$ be the least common
multiple of the denominators of all periods $\{ T_i \}_{i\in[n]}$,
setup times $\{ s_i \}_{i\in[n]}$, and production durations
$\{ \ell_i \}_{i\in[n]}$, written in reduced form. This integer has
polynomial encoding length, since it is bounded by the product of its arguments. We measure schedule lengths in ticks of duration $\frac{ 1 }{ \Delta }$ and write $\hat P=\Delta P$, $\hat T_i=\Delta T_i$, $\hat e_i=\Delta e_i$, and $\hat s_i=\Delta s_i$. These quantities are integers, since $P=\max_{i\in[n]}T_i$ and
$e_i=s_i+\ell_i$. Given a rational query time $t\ge 0$, we may assume without loss of generality that $t \in [0,P)$, noting that it corresponds to forward tick $\lfloor\Delta t\rfloor \in [0, \hat{P}-1]$. Given $\hat P$ ticks numbered from $0$, the latter equals reverse tick $q=\hat P-1-\lfloor\Delta t\rfloor$. Our task is therefore to identify the activity occupying $[q,q+1)$ in reverse time. 

\subsection{Counting free ticks from released work}\label{sec:query-free}

Toward addressing such queries, let us consider an auxiliary machine that receives only a subset $J\subseteq[n]$ of  commodities and works whenever it has unfinished work. We call a tick {$J$-free} if this machine is idle during that tick, and for every integer
$0\le x\le\hat P$, let $F_J(x)$ be the number of $J$-free ticks in $[0,x)$. The important observation is that, to compute $F_J(x)$, we need to know only the amount of work released across times $0,1, \ldots, x-1$, since the order
in which these jobs are processed does not affect $F_J(x)$. In particular, the work released
{strictly before} tick $x$ is
$A_J(x)=\sum_{j\in J}\hat e_j \cdot \lceil\frac{x}{\hat T_j}\rceil$, 
since each commodity $j \in J$ releases jobs of length
$\hat{e}_j$ at times $0,\hat{T}_j, 2\hat{T}_j, \ldots$. The following lemma shows that $F_J(x)$ can be evaluated by knowing the value of $A_J(\cdot)$ only across the set of points $S_J(x)=\{0,x\}\cup
 \{\hat T_j \cdot \lfloor \frac{ x }{ \hat T_j }\rfloor \}_{j\in J}$.

\begin{lemma} \label{lem:query-free-count}
$F_J(x) = \max_{y\in S_J(x)}\{y-A_J(y)\}$, for every integer $0\le x\le\hat P$.
\end{lemma}
\begin{proof}
For any integer $0\le y\le x$, at most $A_J(y)$ ticks of work can be processed by $y$, implying that $F_J(x)\ge F_J(y)\ge y-A_J(y)$, and in turn, $F_J(x) \geq \max_{y\in S_J(x)}\{y-A_J(y)\}$. To prove an inequality in the opposite direction, note that it clearly holds when $F_J(x) = 0$ by choosing $y = 0 \in S_J(x)$ on the right-hand side, and we therefore assume that $F_J(x) \geq 1$. Let $y_\mylast$ be
the right endpoint of the last $J$-free tick in $[0,x)$. By time $y_\mylast$, the
machine has processed all work released strictly before $y_\mylast$, and there
are no $J$-free ticks in $[y_\mylast,x)$. Hence, $F_J(x)=y_\mylast-A_J(y_\mylast)$, and it remains to show that $y_\mylast \in S_J(x)$. 

By definition of $S_J(x)$, the case
$y_\mylast=x$ is immediate, and we therefore assume that $y_\mylast<x$. In this case, at least one job must be released at $y_\mylast$, since tick $y_\mylast$ is busy and the preceding tick was $J$-free. Let $\hat{T}$ be the largest period (measured in ticks) of a commodity
in $J$ releasing a job at $y_\mylast$. Every longer period is a multiple of
$\hat T$, since $\{T_i \}_{i \in [n]}$ are powers of $2$, as explained in Section~\ref{sec:calendar-data}. Consider any commodity $j\in J$ with $\hat T_j>\hat T$.
It cannot release a job at $y_\mylast$, since $\hat T$ was
chosen as the largest period among the commodities releasing
at that time. Moreover, every release time of $j$ is a multiple
of $\hat T_j$, and hence also a multiple of $\hat T$. Since
$y_\mylast$ and $y_\mylast+\hat T$ are consecutive multiples
of $\hat T$, no release of $j$ can occur strictly between them.
Thus, no commodity with period longer than $\hat T$ releases
a job in $[y_\mylast,y_\mylast+\hat T)$. Suppose, for a contradiction, that $y_{\mylast}+\hat T\le x$.
For each commodity $j\in J$ with $\hat T_j\le\hat T$, the period $\hat T_j$
divides $\hat T$, and thus, exactly
$\frac{ \hat T }{ \hat T_j}$ jobs of this commodity are released in $[y_{\mylast},y_{\mylast}+\hat T)$,
beginning with a release at $y_{\mylast}$. The total work released
in this interval is therefore
\[ \sum_{\substack{j\in J:\\\hat T_j\le \hat T}} \frac{\hat T}{\hat T_j} \cdot \hat e_j
 ~~=~~ \hat T \cdot \sum_{\substack{j\in J:\\\hat T_j\le \hat T}}u_j ~~\leq~~ U\hat{T} ~~<~~ \hat T.
\]
Since the preceding tick $[y_\mylast-1,y_\mylast)$ is
$J$-free, all work released strictly before $y_\mylast$ has
already been completed by that time. Consequently, any work
processed in $[y_\mylast,y_\mylast+\hat T)$ must come from
jobs released within this interval, including those released
at its left endpoint. The displayed inequality shows that
their total processing time is strictly less than the
$\hat T$ ticks available. Therefore, at least one tick in
this interval must be $J$-free. Since
$y_\mylast+\hat T\le x$, this observation contradicts the definition of
$y_\mylast$ as the right endpoint of the last $J$-free tick
in $[0,x)$. Thus, $y_{\mylast}\le x<y_{\mylast}+\hat T$, and since $y_{\mylast}$
is a multiple of $\hat T$, it is the last such multiple at or
before $x$. Namely, $y_{\mylast}=\hat T \cdot \lfloor \frac{ x }{ \hat T} \rfloor$, which is one of the points defining $S_J(x)$.
\end{proof}

\subsection{Locating an activity by its priority}\label{sec:query-starts}\label{sec:query-search}

As in Section~\ref{sec:occupation}, we index commodities in
weakly decreasing order of priority, so that $w_1\ge\cdots\ge w_n$,
preserving our tie-breaking order. For convenience, we
introduce a virtual commodity, $n+1$, ranked below all real
commodities, and regard each idle tick as a separate activity of
this commodity, of length $\hat e_{n+1}=1$. For $i=0,\ldots,n$,
let us define the prefix $J_{\leq i}=\{1,\ldots,i\}$ and its complementary suffix
$J_{\geq i+1}=\{i+1,\ldots,n+1\}$. Any {$J_{\geq i+1}$-activity} is
therefore either a block of a real commodity in $J_{\geq i+1}$
or a single idle tick. In the
reverse schedule, jobs of $J_{\leq i}$-commodities may be released
while a $J_{\geq i+1}$-block is in progress, since each block is processed without interruption. However, no such job can be waiting
when a $J_{\geq i+1}$-activity starts, since our placement rule considers all
jobs released by that tick, including those released at the tick
itself, and gives every commodity in $J_{\leq i}$ priority over $J_{\geq i+1}$-commodities. The following claim exploits this observation to
identify the start of each outside activity from the free-tick
count $F_{J_{\leq i}}$.

\begin{lemma}\label{lem:query-free-starts}
Let $a_{\geq i+1}$ be the starting tick of a given $J_{\geq i+1}$-activity, and let $h$ be the total length, in ticks,
of all $J_{\geq i+1}$-activities whose starting ticks lie in
$[0,a_{\geq i+1})$. Then, $F_{ J_{\leq i} }(a_{\geq i+1})=h$ and $F_{ J_{\leq i} }(a_{\geq i+1}+1)=h+1$. 
\end{lemma}
\begin{proof}
Since a $J_{\geq i+1}$-activity starts at $a_{\geq i+1}$,
all $J_{\leq i}$-work released strictly before $a_{\geq i+1}$
has been completed, and no $J_{\leq i}$-job is released at
that tick. The auxiliary machine serving $J_{\leq i}$ has also completed this work, since it receives the
same release sequence and devotes a tick to $J_{\leq i}$ whenever any $J_{\leq i}$-work
remains, whereas the actual machine may be serving another commodity.
Thus, tick $a_{\geq i+1}$ is $J_{\leq i}$-free, and both machines have devoted the same number of ticks to
$J_{\leq i}$ in $[0,a_{\geq i+1})$. Therefore, in this interval, the actual machine spends
$a_{\geq i+1}-h$ ticks on $J_{\leq i}$-work. Since the
auxiliary machine spends the same number of ticks on this work,
it is idle for exactly $h$ ticks in that interval. Hence, $F_{ J_{\leq i} }(a_{\geq i+1})=h$ and $F_{ J_{\leq i} }(a_{\geq i+1}+1)=h+1$. 
\end{proof}

\paragraph{Finding the activity that spans the queried tick.} Let us recall that each block of a real commodity is viewed as a single activity, whereas
each idle tick is a separate activity of the virtual commodity
$n+1$. Since these activities partition $[0,\hat P)$, given a reverse tick
$q\in [0,\hat P-1]$, the activity containing $q$ is the
one whose starting tick is latest among all starts at or before $q$. Algorithm~\ref{alg:query-latest-start} identifies this activity by
visiting the commodities in order of weakly increasing priorities, $n+1,n,\ldots,1$. In each iteration $i$, when considering commodity $i+1$, the commodities visited
by the end of this iteration form the suffix $J_{\geq i+1}$, and we maintain a corresponding record, $(\nu_i,a_i,\delta_i)$. If at least one $J_{\geq i+1}$-activity
has started by $q$, this record specifies the commodity $\nu_i$,
starting tick $a_i$, and length $\delta_i$ in ticks of the latest
such activity to start; otherwise, $(\nu_i,a_i,\delta_i) = (\bot,0,0)$. As explained below, iteration $i$ computes this record from the one at
index $i+1$. In particular, since $J_{\geq 1}$ contains all
commodities, the final record $(\nu_0,a_0,\delta_0)$ identifies
the activity containing $q$.

\medskip

\begin{boxedalgorithm}[label={alg:query-latest-start}]
  {Locating the activity at reverse tick {\boldmath $q$}}
\begin{algorithmic}[1]
\State $(\nu_{n+1},a_{n+1},\delta_{n+1})\gets(\bot,0,0)$
\For{$i=n,n-1,\ldots,0$}
    \State $p\gets F_{J_{\leq i}}(a_{i+1})+\delta_{i+1}$
    \State $m\gets F_{J_{\leq i}}(q+1)$
    \If{$p<m$}
        \State $h\gets p+\hat e_{i+1} \cdot
               \lfloor \frac{ m-1-p }{ \hat e_{i+1} } \rfloor$
        \State Find the smallest $y\in\{0,\ldots,q\}$ with $F_{J_{\leq i}}(y+1)>h$ by binary search
        \State $(\nu_i,a_i,\delta_i)\gets(i+1,y,\hat e_{i+1})$
    \Else
        \State $(\nu_i,a_i,\delta_i)
               \gets(\nu_{i+1},a_{i+1},\delta_{i+1})$
    \EndIf
\EndFor
\State \Return $(\nu_0,a_0,\delta_0)$, identifying commodity $\nu_0$ and reverse interval $[a_0,a_0+\delta_0)$
\end{algorithmic}
\end{boxedalgorithm}

\begin{figure}[htbp!]
\centering
\resizebox{0.95\linewidth}{!}{%
\begin{tikzpicture}[els/base,x=1.2cm,y=1.1cm]


  \foreach \x in {0,...,8} {
    \draw[els/guide]
      (\x,1.08)--(\x,3.10);
    \node[anchor=south,text=elsMid]
      at (\x,3.22) {$\x$};
  }

  \node[anchor=east] at (8.10,3.90) {Reverse time};


  \path[els/block,draw=elsBlue,fill=elsBlue]
    (0,2.40) rectangle (1,3);

  \path[els/block,draw=elsOrange,fill=elsOrange]
    (1,2.40) rectangle (5,3);

  \path[els/block,draw=elsBlue,fill=elsBlue]
    (5,2.40) rectangle (6,3);

  \foreach \x in {6,7}
    \path[els/block,draw=elsMid,fill=elsGray]
      (\x,2.40) rectangle ++(1,.6);

  \node[text=white] at (.5,2.7) {$1$};
  \node[text=white] at (3,2.7) {$2$};
  \node[text=white] at (5.5,2.7) {$1$};
  \node[anchor=east] at (-.20,2.7) {Schedule};

  \fill[elsBlue]
    (3.93,3.16)--(4.07,3.16)--(4,3.04)--cycle;


  \foreach \x in {0,4} {
    \path[els/block,draw=elsBlue,fill=elsBlue]
      (\x,1.2) rectangle ++(1,.6);
    \node[text=white] at (\x+.5,1.5) {$1$};
  }

  \foreach \x in {1,2,3,5}
    \path[els/block,draw=elsOrange,fill=white]
      (\x,1.2) rectangle ++(1,.6);

  \foreach \x in {6,7}
    \path[els/block,draw=elsMid,fill=white]
      (\x,1.2) rectangle ++(1,.6);

  \node[anchor=east] at (-.20,1.5)
    {$J_{\leq 1}=\{1\}$};


  \foreach \x/\h in {1/0,2/1,3/2,5/3,6/4,7/5}
    \draw[els/mapping]
      (\x+.5,1.10)--(\h+1.5,.40);

  \foreach \h in {0,1,2,3} {
    \path[els/block,draw=elsOrange,fill=elsOrange!10]
      (\h+1,-.30) rectangle ++(1,.60);
    \node at (\h+1.5,0) {$\h$};
  }

  \foreach \h in {4,5} {
    \path[els/block,draw=elsMid,fill=elsGray]
      (\h+1,-.30) rectangle ++(1,.60);
    \node at (\h+1.5,0) {$\h$};
  }

  \node[anchor=east,align=right] at (.78,0)
    {$J_{\leq 1}$-free-tick\\index};

  \draw[elsOrange,line width=.55pt]
    (1,-.60)--(5,-.60);
  \draw[elsOrange,line width=.55pt]
    (1,-.72)--(1,-.48);
  \draw[elsOrange,line width=.55pt]
    (5,-.72)--(5,-.48);

  \node[anchor=north,text=elsOrange] at (3,-.73)
    {$\hat e_2=4$};

\end{tikzpicture}%
}
\caption{The reverse schedule and the auxiliary machine serving
only commodity~1, which has higher priority. Here,
$(\hat T_1,\hat T_2)=(4,8)$ and $(\hat e_1,\hat e_2)=(1,4)$.
Gray cells in the actual schedule are idle ticks, and white cells
in the auxiliary schedule are $J_{\leq 1}$-free ticks. The marker at tick~4
denotes a new commodity-1 release. The four indices assigned to
block $[1,5)$ correspond to free ticks $1$, $2$, $3$, and $5$, which need not be consecutive in time.}
\label{fig:free-tick-correspondence}
\end{figure}
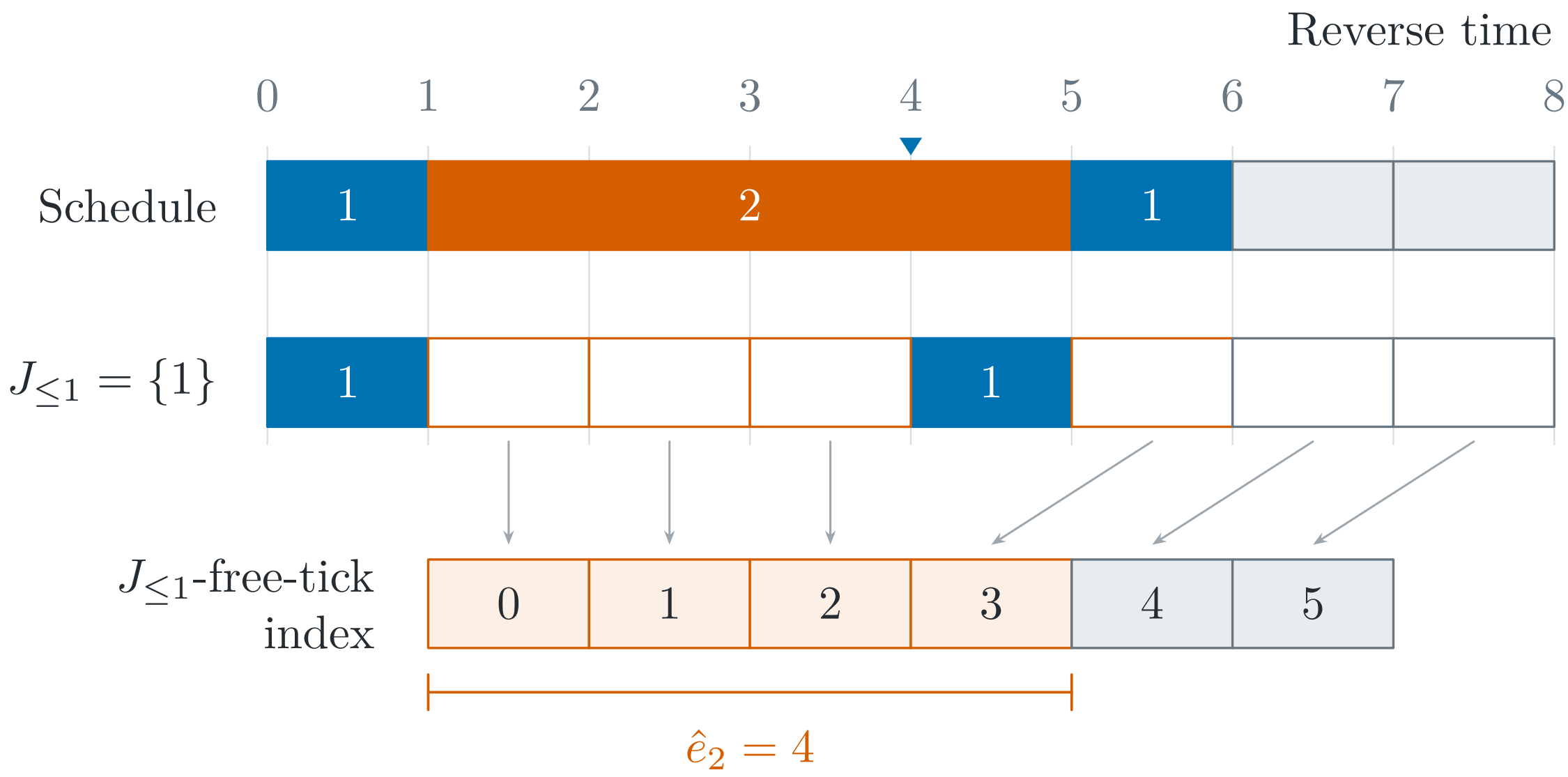

\paragraph{Interpreting our update procedure.} To better understand the inner workings of our algorithm, let us focus on a single iteration $i$, and number $J_{\leq i}$-free ticks from
$0$. By Lemma~\ref{lem:query-free-starts}, if a
$J_{\geq i+1}$-activity starts at tick $a$, then $a$ is the $J_{\leq i}$-free
tick with index $F_{J_{\leq i}}(a)$. If this activity has length $d$, then the total length of all
$J_{\geq i+1}$-activities up to and including it is
$F_{J_{\leq i}}(a)+d$. By the same lemma, this is the index of
the $J_{\leq i}$-free tick at which the next $J_{\geq i+1}$-activity starts,
if there is one. We therefore associate the current activity
with the $d$ consecutive indices $F_{J_{\leq i}}(a),\ldots,F_{J_{\leq i}}(a)+d-1$. The first index identifies the
actual starting tick of the activity; the free ticks in its group
need not be consecutive in the schedule. The first group starts at index $0$, since no
$J_{\geq i+1}$-activity precedes the first one, and each
subsequent group starts immediately after the preceding group.
Moreover, the actual machine and the auxiliary machine complete
the same amount of $J_{\leq i}$-work over the entire period. Hence, the
total length of $J_{\geq i+1}$-activities equals the number
of $J_{\leq i}$-free ticks, so the groups exhaust the list. Figure~\ref{fig:free-tick-correspondence} illustrates this grouping scheme for $J_{\leq 1}=\{1\}$. The block of commodity~2 starts at the
first free tick assigned to its group, but the remaining free ticks
in that group need not occur during the block.

We proceed by elaborating on the role of $p$ and $m$ in iteration $i$. If $\nu_{i+1}\ne\bot$, the previously
recorded activity starts at $a_{i+1}$ and has length $\delta_{i+1}$.
Therefore, the first index after its group in the current list is 
$p=F_{J_{\leq i}}(a_{i+1})+\delta_{i+1}$. If $\nu_{i+1}=\bot$, there is no recorded group to skip,
$a_{i+1}=\delta_{i+1}=0$, and the same formula leads to $p=0$. In turn, $m=F_{J_{\leq i}}(q+1)$ stands for the number of $J_{\leq i}$-free ticks through $q$. As such, a group has started by $q$ precisely when its first
index is less than $m$, and it remains to find the last group, if any,
whose first index lies between $p$ and $m-1$.

\paragraph{Why the priority scan works.}
We show that each iteration produces the stated record. Before
the first iteration, no commodity has been visited, and the empty
record at index $n+1$ is correct. Now, suppose that
$(\nu_{i+1},a_{i+1},\delta_{i+1})$ is correct for the commodities
visited before iteration $i$. If $p\ge m$, no group with first index at least $p$ has started
by $q$. Thus, if the previous record is nonempty, no later
$J_{\geq i+1}$-activity has started by $q$; if it is empty,
no $J_{\geq i+1}$-activity has started by $q$ at all. In either
case, keeping the previous record is correct.

Hence, the non-trivial scenario to consider is $p<m$. Here, the group beginning at index $p$ and
every subsequent group that has started by $q$ must belong to
commodity $i+1$. Indeed, all other commodities in $J_{\geq i+1}$
have already been visited. If a group of any previously visited commodity had its first
index between $p$ and $m-1$, its activity would have started
by $q$. For a nonempty previous record, this start would be
later than $a_{i+1}$, contrary to the choice of the record.
For an empty previous record, it would contradict the fact
that no previously visited commodity had started an activity
by $q$. Each of these groups has length $\hat e_{i+1}$, meaning that their first
indices are $p,p+\hat e_{i+1},p+2\hat e_{i+1},\ldots$, up to the
last such index below $m$, which is exactly $h= p+\hat e_{i+1} \cdot \lfloor \frac{ m-1-p }{ \hat e_{i+1} } \rfloor$ in
Algorithm~\ref{alg:query-latest-start}. By
Lemma~\ref{lem:query-free-starts}, the corresponding activity
starts at the smallest integer $y$ for which
$F_{J_{\leq i}}(y+1)>h$. Since this activity has started by $q$,
the binary search recovers its starting tick in $\{0,\ldots,q\}$.
Thus, $(\nu_i,a_i,\delta_i)=(i+1,y,\hat e_{i+1})$ is the correct
record for $J_{\geq i+1}$. After iteration $0$, the suffix $J_{\geq 1}$ contains
all commodities, including the virtual commodity, and
$[a_0,a_0+\delta_0)$ is the reverse interval spanning tick $q$.

\paragraph{Running time.}
To obtain the final record $(\nu_0,a_0,\delta_0)$,
Algorithm~\ref{alg:query-latest-start} performs $n+1$ iterations.
In each iteration $i$, it evaluates $F_{J_{\leq i}}(a_{i+1})$ and
$F_{J_{\leq i}}(q+1)$ to compute $p$ and $m$. If $p<m$, it also performs a binary
search for the new starting tick, involving $O(\log{\hat P})$ further
evaluations. By Lemma~\ref{lem:query-free-count}, each evaluation
of $F_J(x)$ takes $O(n^2)$ arithmetic operations, since there are at most
$|J|+2$ candidate points in $S_J(x)$, and computing $A_J(y)$ at
each point requires summing $|J|$ terms. If the final record
identifies an idle tick, two additional binary searches recover the
maximal idle interval containing $q$. Its endpoints are the first
and last integer boundaries $x\in[0,\hat P]$ at which the busy-tick
count $x-F_{[n]}(x)$ equals $q-F_{[n]}(q)$. Converting the resulting
reverse interval to absolute time and determining whether $t$ falls
in setup or production for a real commodity takes $O(1)$ additional
operations. Thus, an activity query requires
$O(n^3\log{\hat P})$ arithmetic operations. Finally, all tick
boundaries and released-work counts are at most $\hat P$, and 
$\hat P=P\Delta$ has polynomial encoding length, as argued in Sections~\ref{sec:calendar-data} and~\ref{sec:query-reverse}. Consequently, the entire query
takes time polynomial in the binary encoding length of the given instance, as claimed in Theorem~\ref{thm:activity-queries}.
\section{Concluding Remarks}\label{sec:conclusion}

We conclude this paper by highlighting a number of fundamental open questions for future research. Specifically, we first describe prospective refinements of our construction, potentially saving a large portion of its approximation guarantee. We then ask about the exact strength of our convex relaxation, the regularity required for answering activity queries, and conditions under which shorter cycles could suffice.

\paragraph{Improved approximation guarantees.} As previously mentioned, we intentionally do not optimize any constants, focusing instead on simplicity and
accessibility. Interestingly, our subsequent investigation suggests a
possible route to an improved approximation ratio of roughly
$8.48$. We regard these initial attempts as a promising direction, although the
emerging construction and its analysis appear substantially more
involved, and it is still unclear whether they can support polynomial-time
activity queries. We therefore prefer
to concentrate on explaining how competition for machine time creates inventory
requirements, captured within a convex program that subsequently guides the construction
of a provably good production policy, coupled with a compact description from
which its individual activities can be efficiently recovered. Establishing clean and digestible improvements on our work is a compelling open question for future research.

\paragraph{The gap between $\OPT\eqref{eq:inventory-convex}$ and $\OPTdyn$.} Our approximation factor combines two intermixed losses: The gap between the convex program~\eqref{eq:inventory-convex} and the dynamic optimum, and the extra cost incurred when converting a relaxation solution into a cyclic policy. While this relaxation records inventory carried during other commodities' blocks and enforces consistency through pairwise inequalities, it does not specify a common chronological order in which  these blocks occur. Does this missing information produce a substantial gap, or do the existing constraints already capture most of the inventory cost imposed by sharing a common machine? It would be interesting to either construct instances exhibiting a non-trivial gap or prove tighter upper bounds. One could also incorporate stronger group inequalities that are not implied by the present pairwise constraints, provided that the extended formulation remains efficiently solvable. 

\paragraph{How much regularity do activity queries need?} A close inspection of the  start-time identity in Lemma~\ref{lem:query-free-starts} reveals that it applies to arbitrary integral periods, due to identifying the start of an activity via the free-tick counts of an auxiliary machine serving higher-priority commodities. Here, the arithmetic period structure enters when free-tick counts are computed and inverted to recover a particular tick, with powers of $2$ making these operations relatively simple. Does an equally efficient procedure exist for unrestricted periods that do not form a divisible chain, possibly with many distinct divisors? It would also be useful to understand whether comparably compact representations are possible with shifted release times or alternative scheduling rules that offer better holding-cost guarantees.

\paragraph{When do shorter cycles suffice?}
The two-commodity construction in Appendix~\ref{sec:cycle-length} shows that neither a fixed number of commodities nor utilization bounded away from $1$ ensures non-trivial cyclic policies with polynomially many runs. Its parameters, however, have widely separated demand loads and setup costs. This observation suggests asking whether restrictions on the numerical ranges of the underlying parameters imply substantially shorter approximate cycles. For example, what can be guaranteed when the ratios between the largest and smallest values of $\rho_i$, $K_i$, and $c_i$ are polynomially bounded in the input length? A related question asks about representations falling between explicit lists and activity queries: Can structured policies with comparable approximation guarantees be described by polynomially many run-and-idle patterns and repetition counts written in binary? The specific instance in Appendix~\ref{sec:cycle-length} clearly admits this type of description; however, we still do not know whether it is possible in general.

\paragraph{Acknowledgment of AI assistance.}
ChatGPT~6.0 Astra assisted in preparing Figures~\ref{fig:production-inventory-dynamics}--\ref{fig:free-tick-correspondence},
identifying relevant literature, and refining 
language and exposition. The author is fully responsible
for the final manuscript, including the accuracy of its
mathematical results, references, and figures.

\phantomsection
\addcontentsline{toc}{section}{References}
\bibliographystyle{plainnat}
\bibliography{BIB-ELSP}

@article{AdelmanBarz2014,
  author  = {Adelman, Daniel and Barz, Christiane},
  title   = {A Unifying Approximate Dynamic Programming Model for the Economic Lot Scheduling Problem},
  journal = {Mathematics of Operations Research},
  year    = {2014},
  volume  = {39},
  number  = {2},
  pages   = {374--402}
}

@article{NobilEtAl2026,
  author  = {Nobil, Amir Hossein and
             Afshar Sedigh, Amir Hosein and
             C{\'a}rdenas-Barr{\'o}n, Leopoldo Eduardo and
             Nobil, Erfan and
             C{\'e}spedes-Mota, Armando and
             Trevi{\~n}o-Garza, Gerardo},
  title   = {A review paper on the economic lot scheduling problem
             with the common cycle approach},
  journal = {RAIRO -- Operations Research},
  year    = {2026},
  volume  = {60},
  number  = {3},
  pages   = {1325--1355}
}

@article{DollWhybark1973,
  author  = {Doll, Charles L. and Whybark, David C.},
  title   = {An Iterative Procedure for the Single-Machine
             Multi-Product Lot Scheduling Problem},
  journal = {Management Science},
  year    = {1973},
  volume  = {20},
  number  = {1},
  pages   = {50--55}
}

@book{NahmiasOlsen2015,
  author    = {Nahmias, Steven and Olsen, Tava Lennon},
  title     = {Production and Operations Analysis},
  edition   = {seventh},
  publisher = {Waveland Press},
  address   = {Long Grove, Illinois},
  year      = {2015}
}

@article{Hsu1983,
  author  = {Hsu, Wen-Lian},
  title   = {On the General Feasibility Test of Scheduling Lot Sizes
             for Several Products on One Machine},
  journal = {Management Science},
  year    = {1983},
  volume  = {29},
  number  = {1},
  pages   = {93--105}
}

@article{Bomberger1966,
  author  = {Bomberger, Earl E.},
  title   = {A Dynamic Programming Approach to a Lot Size Scheduling Problem},
  journal = {Management Science},
  year    = {1966},
  volume  = {12},
  number  = {11},
  pages   = {778--784}
}

@article{Carstensen1999,
  author  = {Carstensen, Peter},
  title   = {The Economic Lot Scheduling Problem --- Survey and {LP}-Based Method},
  journal = {OR Spektrum},
  year    = {1999},
  volume  = {21},
  number  = {4},
  pages   = {429--460}}

@article{Dobson1987,
  author  = {Dobson, Gregory},
  title   = {The Economic Lot-Scheduling Problem: Achieving Feasibility Using Time-Varying Lot Sizes},
  journal = {Operations Research},
  year    = {1987},
  volume  = {35},
  number  = {5},
  pages   = {764--771}
}

@article{Elmaghraby1978,
  author  = {Elmaghraby, Salah E.},
  title   = {The Economic Lot Scheduling Problem ({ELSP}): Review and Extensions},
  journal = {Management Science},
  year    = {1978},
  volume  = {24},
  number  = {6},
  pages   = {587--598}
}

@article{GallegoShaw1997,
  author  = {Gallego, Guillermo and Shaw, Dong Xiao},
  title   = {Complexity of the {ELSP} with General Cyclic Schedules},
  journal = {IIE Transactions},
  year    = {1997},
  volume  = {29},
  number  = {2},
  pages   = {109--113}
}

@incollection{Lovasz2003SDP,
  author    = {Lov{\'a}sz, L{\'a}szl{\'o}},
  title     = {Semidefinite Programs and Combinatorial Optimization},
  booktitle = {Recent Advances in Algorithms and Combinatorics},
  editor    = {Reed, Bruce A. and Sales, Cl{\'a}udia L.},
  year      = {2003},
  pages     = {137--194},
  publisher = {Springer},
  address   = {New York}
}

@article{Roundy1989,
  author  = {Roundy, Robin},
  title   = {Rounding Off to Powers of Two in Continuous Relaxations of Capacitated Lot Sizing Problems},
  journal = {Management Science},
  year    = {1989},
  volume  = {35},
  number  = {12},
  pages   = {1433--1442}
}

@article{Zipkin1991,
  author  = {Zipkin, Paul H.},
  title   = {Computing Optimal Lot Sizes in the Economic Lot Scheduling Problem},
  journal = {Operations Research},
  year    = {1991},
  volume  = {39},
  number  = {1},
  pages   = {56--63}
}

@book{SilverPykeThomas2017,
  author    = {Silver, Edward A. and Pyke, David F. and Thomas, Douglas J.},
  title     = {Inventory and Production Management in Supply Chains},
  edition   = {fourth},
  year      = {2017},
  publisher = {CRC Press},
  address   = {Boca Raton, FL}
}

@article{JonesInman1989,
  author  = {Jones, Philip C. and Inman, Robert R.},
  title   = {When Is the Economic Lot Scheduling Problem Easy?},
  journal = {IIE Transactions},
  year    = {1989},
  volume  = {21},
  number  = {1},
  pages   = {11--20}
}

@techreport{Jackson1955,
  author      = {Jackson, James R.},
  title       = {Scheduling a Production Line to Minimize Maximum Tardiness},
  institution = {Management Science Research Project, University of California, Los Angeles},
  type        = {Research Report},
  number      = {43},
  year        = {1955}
}

@article{GrahamEtAl1979,
  author  = {Graham, Ronald L. and Lawler, Eugene L. and Lenstra, Jan K. and Rinnooy Kan, Alexander H. G.},
  title   = {Optimization and Approximation in Deterministic Sequencing and Scheduling: A Survey},
  journal = {Annals of Discrete Mathematics},
  year    = {1979},
  volume  = {5},
  pages   = {287--326}
}

@article{Rogers1958,
  author  = {Rogers, Jack},
  title   = {A Computational Approach to the Economic Lot Scheduling Problem},
  journal = {Management Science},
  year    = {1958},
  volume  = {4},
  number  = {3},
  pages   = {264--291}
}

@article{Maxwell1964,
  author  = {Maxwell, William L.},
  title   = {The Scheduling of Economic Lot Sizes},
  journal = {Naval Research Logistics Quarterly},
  year    = {1964},
  volume  = {11},
  number={2},
  pages   = {89--124}
}

@article{ChanEtAl2013,
  author  = {Chan, Hing Kai and Chung, Sai Ho and Lim, Ming K.},
  title   = {Recent Research Trend of Economic-Lot Scheduling Problems},
  journal = {Journal of Manufacturing Technology Management},
  year    = {2013},
  volume  = {24},
  number  = {3},
  pages   = {465--482}
}

@article{SantanderMercadoJubizDiaz2016,
  author  = {Santander-Mercado, Alcides and Jubiz-Diaz, Maria},
  title   = {The Economic Lot Scheduling Problem: A Survey},
  journal = {International Journal of Production Research},
  year    = {2016},
  volume  = {54},
  number  = {16},
  pages   = {4973--4992}
}

@article{BeckGlock2020,
  author  = {Beck, Fabian G. and Glock, Christoph H.},
  title   = {The Economic Lot Scheduling Problem: A Content Analysis},
  journal = {International Journal of Production Research},
  year    = {2020},
  volume  = {58},
  number  = {11},
  pages   = {3437--3454}
}

@book{Zipkin2000,
  author    = {Zipkin, Paul H.},
  title     = {Foundations of Inventory Management},
  year      = {2000},
  publisher = {McGraw-Hill},
  address   = {Boston, MA}
}

@book{Pinedo2005,
  author    = {Pinedo, Michael L.},
  title     = {Planning and Scheduling in Manufacturing and Services},
  year      = {2005},
  publisher = {Springer},
  address   = {New York, NY}
}

@book{Axsater2015,
  author    = {Axs{\"a}ter, Sven},
  title     = {Inventory Control},
  edition   = {Third},
  year      = {2015},
  publisher = {Springer}
}

\addtocontents{toc}{\protect\setcounter{tocdepth}{2}} 
\appendix
\section{Additional Proofs from Section~\ref{sec:intro}}

\subsection{Feasibility condition} \label{app:feasibility}

\paragraph{Sufficiency.} Suppose either $R<1$, or $R=1$ and $s_1 = \cdots = s_n =0$. In this scenario, we explicitly define a cyclic policy $\pi$ by first choosing some period $P>0$ for which $\kappa P\ge\sum_{i\in[n]}s_i$. In each such period, we visit the commodities in a fixed order, setting up each commodity $i \in [n]$ once and producing it for $\rho_iP$ units of time. These operations occupy  $\sum_{i\in[n]}s_i+RP\le P$ units of time, with the remainder left idle. Now, with initial inventory $I_i^\pi(0)=d_iP$, it is easy to verify that $\pi$ is a feasible policy. Moreover, since the run of commodity $i$ produces $p_i \rho_i P = d_i P$ units, every inventory returns to its initial level at time $P$, implying that $\pi$ has finite long-run average cost.

\paragraph{Necessity.}
Suppose that $\pi$ is a feasible policy. When each commodity's
inventory equation~\eqref{eqn:inventory_time} is divided by its production rate, summing over all commodities leads to
\begin{equation} \label{eqn:proof_feasibility_1}
\sum_{i\in[n]}\frac{I_i^\pi(t)}{p_i}
 ~~=~~ \sum_{i\in[n]} \left( \frac{I_i^\pi(0)}{p_i}  +|\mathcal P_i\cap[0,t)|- \frac{ d_i t }{ p_i } \right) ~~\le~~  \sum_{i\in[n]}\frac{I_i^\pi(0)}{p_i}+(1-R)t. 
\end{equation}
Here, the inequality holds since the production-time sets $\{ {\cal P}_i \}_{i \in [n]}$ are pairwise disjoint, implying that $\sum_{i\in[n]}   |\mathcal P_i\cap[0,t)| \leq t$. If $R>1$, the right-hand side becomes negative
for sufficiently large $t$, since the initial inventories $(I_i^\pi(0))_{i \in [n]}$ are
finite, which would contradict the feasibility of $\pi$.

Knowing that $R\le1$, it remains to show that $R=1$ requires all setup times to be
zero. For this purpose, note that 
\[
 0
 ~~\le~~
 t-\sum_{i\in[n]}|\mathcal P_i\cap[0,t)|
 ~~=~~
 \sum_{i\in[n]}\frac{I_i^\pi(0)-I_i^\pi(t)}{p_i}
 ~~\le~~
 \sum_{i\in[n]}\frac{I_i^\pi(0)}{p_i},
\]
where the sole equality above is exactly the one in~\eqref{eqn:proof_feasibility_1} with $R=1$. As a result, the total time $ t-\sum_{i\in[n]}|\mathcal P_i\cap[0,t)|$ spent on setups
and idling during $[0,t)$ is bounded by a fixed
finite quantity, uniformly over all $t\ge0$. Now, suppose there exists a commodity $i\in[n]$ with setup time $s_i > 0$. Then, due to requiring infinitely many production
runs, $\pi$ has infinitely many nonoverlapping setups of this commodity, each of duration
$s_i > 0$. Therefore, their cumulative duration cannot be bounded,
contradicting the fixed upper bound on $ t-\sum_{i\in[n]}|\mathcal P_i\cap[0,t)|$.

\subsection{Number of runs: Exponential lower bound} \label{sec:cycle-length}

\paragraph{Instance description.} We construct a family of two-commodity instances indexed by an integer
$k\ge1$. Both commodities have unit production and holding rates, i.e., $p_1 = p_2 = 1$ and $h_1 = h_2 = 1$, as well as identical setup times $s_1 = s_2 =
\frac{ 1 }{ 16 \cdot 4^k}$. Their demand rates and setup costs are $(d_1,K_1)=(\frac{ 1 }{ 4 }, \frac{ 1 }{ 4^k})$ and $(d_2,K_2)=(\frac{ 1 }{ 4 \cdot 4^k},1)$. 
All data are positive and have total binary encoding length $\Theta(k)$. In this setting, we show that for any $\alpha \geq 1$, every feasible cyclic policy $\pi$ with long-run average cost $\mathcal C(\pi) \le \alpha \cdot \OPTdyn$ must have at least $\frac{4^k}{24\alpha^2}$ production runs per cycle.

\paragraph{Upper-bounding $\OPTdyn$.} Our first claim is that $\OPTdyn \leq \frac{ 3 }{ 2^k }$. For this purpose, we describe a cyclic policy $\pi_k$ with long-run average cost $\mathcal{C}( \pi_k ) \leq \frac{ 3 }{ 2^k }$. This policy divides a period of length $2^k$ into $4^k$ slots of length $\frac{ 1 }{ 2^k }$. In each slot, commodity~1 is set up and produced for $\frac{ 1 }{ 4 \cdot 2^k }$ units of time. In the first slot only, commodity~2 has a similar setup and run length, which is possible since $2 \cdot ( \frac{ 1 }{ 16 \cdot 4^k} + \frac{ 1 }{ 4 \cdot 2^k } ) < \frac{ 1 }{ 2^k }$. These choices ensure that each commodity's output balances its demand between successive runs. We set the initial stock of each commodity so that inventory reaches zero when its production begins. As a result, its stock stays between zero and the entire output of a single run, $\frac{ 1 }{ 4 \cdot 2^k }$. Each commodity incurs a setup cost of $\frac{ 1 }{ 2^k }$ per unit time, and since the holding rates are $h_1 = h_2 = 1$, even charging each commodity for its maximum inventory throughout the cycle leads to a holding cost of at most $\frac{ 1 }{ 4 \cdot 2^k }$ per unit time, meaning that $\mathcal{C}( \pi_k ) \leq 2 \cdot ( \frac{ 1 }{ 2^k } + \frac{ 1 }{ 4 \cdot 2^k }) < \frac{ 3 }{ 2^k }$.

\paragraph{The number of runs in an approximate cycle.} Now, consider an $\alpha$-approximate cyclic policy $\pi$ of period $P$ with $m$ runs in total. Commodity~2 must have at least one run per period and therefore costs at least $\frac{ 1 }{ P }$ in setups per unit time. Since $d_1 = \frac{ 1 }{ 4 }$ and $p_1 = 1$, commodity~1 must occupy $\frac{ P }{ 4 }$ units of production time in at most $m$ runs, meaning that its mean run length is at least $\frac{ P }{ 4m }$. As such, Lemma~\ref{lem:individual} bounds its mean holding cost from below by $(1-\frac{ 1 }{ 4 }) \cdot \frac{ P }{ 8m }= \frac{ 3P }{ 32m }$. By adding these two costs, ${\cal C}(\pi) \geq \frac1P+\frac{3P}{32m}\ge\sqrt{\frac3{8m}}$, via the AM--GM inequality. On the other hand, given our upper bound on $\OPTdyn$, we have $\mathcal C(\pi) \leq \alpha \cdot \OPTdyn \leq \frac{ 3\alpha }{ 2^k }$, and therefore, $m \geq \frac{4^k}{24\alpha^2}$.

\subsection{The classical relaxation: Unbounded gap}\label{sec:gap}

\paragraph{Instance description.} To demonstrate an unbounded gap for the classical frequency relaxation~\eqref{eq:bomberger}, let us consider the following two-commodity instance. Both commodities have production rate $1$ and demand rate $\frac{ 1 }{ 4 }$, meaning that $\kappa = 1-(\rho_1 + \rho_2)=\frac{ 1 }{ 2 }$. Their setup times are $\delta^2$, for a rational parameter $0<\delta\le \frac{ 1 }{ 2 }$. Finally, we fix setup costs and holding rates to $(K_1,h_1)=(\delta^2,1)$ and $(K_2,h_2)=(1,\delta^2)$. Thus, commodity~1 is inexpensive to set up but expensive to hold,
whereas commodity~2 has the opposite cost structure. In what follows, we argue that $L_{\rm B} \leq \frac{5\delta}{4}$ whereas $\OPTdyn \ge \frac{1}{4\sqrt2}$, leading to a ratio of $\frac{\OPTdyn}{L_{\rm B}} \ge \frac{1}{5\sqrt2\delta}$, which   tends to infinity as $\delta$ tends to $0$.

\paragraph{Upper bound on $L_{\rm B}$.} Consider the candidate solution $(f_1, f_2) = (\frac{1}{4\delta}, \frac{\delta}{4})$. These frequencies satisfy the average setup-time constraint $\sum_{i\in[n]} s_i f_i \le \kappa$, since $s_1f_1+s_2f_2 = \frac{\delta+\delta^3}{4} \le \frac12
= \kappa$. Therefore, $(f_1, f_2)$ is a feasible solution to relaxation~\eqref{eq:bomberger}. Moreover, both commodities have setup terms $K_if_i = \frac{ \delta }{ 4 }$ and holding terms $\frac{h_id_i(1-\rho_i)}{2f_i} = \frac{ 3\delta }{ 8 }$, leading to $L_{\rm B} \leq 2 \cdot (\frac{\delta}{4}+\frac{3\delta}{8})  = \frac{5\delta}{4}$.

\paragraph{Lower bound on $\OPTdyn$.}
To derive this lower bound, let us first consider any feasible cyclic policy $\pi$ of finite cost, with period length $P$. Let $m\ge1$ be its number of production runs for
commodity~2, and let $\ell_1,\ldots,\ell_m$ be their durations, excluding setups. Production must balance demand
over a single period, leading to
$\sum_{r \in [m]} \ell_r= \frac{ P }{ 4 }$. During a run of commodity~2 of length $\ell$, commodity~1 is not produced and its inventory decreases at rate $\frac{ 1 }{ 4}$.
Since its inventory must be nonnegative at the end of this run,
at time $u\in[0,\ell]$ after the run begins, it must be at least
$\frac{ \ell-u }{ 4 }$. Since $h_1=1$, the holding cost incurred by commodity~1
during this run alone is at least $\int_0^\ell \frac{\ell-u}{4}\dd u = \frac{\ell^2}{8}$. By taking into account only these holding costs and the setups of commodity~2, we infer that
\[ \mathcal C(\pi) ~~\ge~~ \frac{m}{P}   +\frac{1}{8P} \cdot \sum_{r \in [m]}\ell_r^2 ~~\ge~~ \frac{m}{P} +\frac{1}{8Pm} \cdot \left( \sum_{r \in [m]}\ell_r \right)^2 ~~=~~ \frac{m}{P}+\frac{P}{128m} ~~\ge~~ \frac{1}{4\sqrt2}. \]
Here, the second inequality follows from Cauchy--Schwarz, the equality substitutes $\sum_{r \in [m]} \ell_r= \frac{ P }{ 4 }$, and the last inequality follows from the AM--GM inequality. This lower bound holds for every feasible cyclic policy of finite cost. By Lemma~\ref{lem:cyclic-density}, the infimum over such policies equals the unrestricted dynamic optimum, and hence
$\OPTdyn \ge \frac{1}{4\sqrt2}$.

\section{Additional Proofs from Section~\ref{sec:lower}}

\subsection{Proof of Lemma~\ref{lem:cyclic-density}}
\label{sec:cyclic-density-proof}

Since every cyclic policy is also dynamic, we have $\OPTdyn\le\OPTcyc$. To prove the reverse inequality, let us 
fix a feasible dynamic policy $\pi\in\Pi$ of finite long-run
average cost. For every $\eps>0$, we show that there exists a feasible
cyclic policy $\hat{\pi}\in\Pi_{\mathrm{cyc}}$ for which $ \mathcal C(\hat{\pi})\le\mathcal C(\pi)+\eps$. 

\paragraph{Segment with nearby endpoint inventories.}
Letting $C=\mathcal C(\pi)$, there is an increasing sequence
$(t_k)_{k\in\bbN}$ tending to infinity such that
$\sum_{i\in[n]}h_i I_i^\pi(t_k)\le C+1$ for every $k\in\bbN$.
Otherwise, there would be some $T_0 \geq 0$ such that $\pi$ has a holding cost rate of $\sum_{i\in[n]}h_i I_i^\pi(t)\geq C+1$ for all $t > T_0$, forcing its long-run average cost to be at least $C+1$. We next adjust this sequence so that none of its elements lies in the interior of a setup-production block. If $t_k$ lies in the interior of such a block, say of commodity $i$, let $\tau_k$ be the beginning of that block; otherwise,
$\tau_k=t_k$. As a result, if $t_k$ occurs during setup, it moves back by
at most $s_i$. Alternatively, if production has been underway for $u$ time units,
then $I_i^\pi(t_k)\ge(p_i-d_i)u$, since the inventory of this commodity was
nonnegative when production began. As $h_iI_i^\pi(t_k)\le C+1$,
we obtain the uniform bound $0 \le t_k-\tau_k \le \max_{i\in[n]} \{s_i+\frac{C+1}{h_i(p_i-d_i)}\}$, implying in particular that $(\tau_k)_{k \in \bbN}$ tends to infinity. Moreover, since the inventory of each
commodity $i$ changes at a rate of at most
$\max\{d_i,p_i-d_i\}$ in absolute value, moving from $t_k$
to $\tau_k$ changes each inventory by a bounded amount, and the
vectors $(I^\pi(\tau_k))_{k \in \bbN}$ remain uniformly bounded.

Since $(\tau_k)_{k \in \bbN}$ tends to infinity, there exists a subsequence on
which the times $\tau_k$ are strictly increasing. Along this
subsequence, the inventory vectors $I^\pi(\tau_k)$ lie in a
fixed compact box in $\mathbb R^n$, since each commodity's
inventory is uniformly bounded. By the Bolzano--Weierstrass
theorem, these vectors admit a convergent subsequence, whose
times are still strictly increasing. Hence,
for any $\delta>0$, we can fix one endpoint $a$ sufficiently far
along this subsequence so that every subsequent endpoint $b$
satisfies $|I_i^\pi(b)-I_i^\pi(a)|\le\delta$ for all $i\in[n]$.

Let $b_1,b_2,\ldots$ be the successive endpoints after $a$ in
the chosen subsequence, so that $(b_\ell-a)_{\ell \in \bbN}$ tends to infinity. To upper-bound the average cost $\pi$ incurs over segments of the form $[a, b_{\ell})$, let ${\cal C}_{[0,T)}(\pi)$ denote the total cost
incurred across $[0,T)$, namely,
\[
 {\cal C}_{[0,T)}(\pi)
 ~~=~~ \sum_{i\in[n]}K_i N_i^\pi(T)
   +\int_0^T\sum_{i\in[n]}h_i I_i^\pi(t)\dd t.
\]
In addition, for $0\le T_1 \le T_2$, we write
${\cal C}_{[T_1,T_2)}(\pi)={\cal C}_{[0,T_2)}(\pi)-{\cal C}_{[0,T_1)}(\pi)$. Then, 
\begin{align*}
 \limsup_{\ell\to\infty}
 \frac{{\cal C}_{[a,b_\ell)}(\pi)}{b_\ell-a}
 ~~&=~~ \limsup_{\ell\to\infty}
    \frac{{\cal C}_{[0,b_\ell)}(\pi)-{\cal C}_{[0,a)}(\pi)}
         {b_\ell-a} \\
 &=~~ \limsup_{\ell\to\infty}
    \frac{{\cal C}_{[0,b_\ell)}(\pi)}{b_\ell}
 ~~\le~~ \limsup_{T\to\infty}
    \frac{{\cal C}_{[0,T)}(\pi)}{T}
 ~~=~~ C .
\end{align*}

\paragraph{Balancing production over the cycle.} Let us fix $\tau>0$, which is a parameter whose value will be determined later on, and choose
$\delta>0$ small enough that $\delta \le \tau \cdot \min_{i\in[n]}d_i$ and $\delta \cdot \sum_{i\in[n]}\frac{1}{p_i}\le\kappa\tau$.
For this choice of $\delta$, let $a$ and $(b_\ell)_{\ell \in \bbN}$
be the endpoints obtained in the previous paragraph. For each
sufficiently large $\ell$, we will extend the segment $[a,b_\ell)$
to a period of length $P_\ell=b_\ell-a+\tau$ over which production
balances demand for every commodity.

By inventory balance, across $[a,b_\ell)$, the production quantity of commodity $i$ is $d_i(b_\ell-a)+I_i^\pi(b_\ell)-I_i^\pi(a)$.
Since $(b_\ell-a)_{\ell \in \bbN}$ tends to infinity and
$|I_i^\pi(b_\ell)-I_i^\pi(a)|\le\delta$, this quantity is strictly positive
for every $i\in[n]$ once $\ell$ is sufficiently large. By construction, neither
endpoint splits a block, meaning that the segment $[a,b_\ell)$ contains complete runs. To bring the production quantity of commodity $i$ to
$d_iP_\ell$, we extend its first production run by
$x_{i\ell} =\frac{d_i\tau+I_i^\pi(a)-I_i^\pi(b_\ell)}{p_i}$
units of time, which is nonnegative since
$I_i^\pi(b_\ell)-I_i^\pi(a)\le\delta\le d_i\tau$. Moreover,
all extensions combined amount to at most $\tau$ units of time, since
\[
 \tau-\sum_{i\in[n]}x_{i\ell}
 ~~=~~(1-R)\tau+
   \sum_{i\in[n]}\frac{I_i^\pi(b_\ell)-I_i^\pi(a)}{p_i}
 ~~\ge~~\kappa\tau- \delta \cdot \sum_{i\in[n]}\frac{1}{p_i}
 ~~\ge~~0.
\]
We insert the remaining $\tau-\sum_{i\in[n]}x_{i\ell}$ time as an idle interval at the beginning
of the segment, keeping all setups and runs in their original
order. The resulting schedule has length $P_\ell = b_{\ell} - a + \tau$ and requires
no additional setups. Its total output of commodity $i$ is $d_i(b_\ell-a)+I_i^\pi(b_\ell)-I_i^\pi(a)+p_i x_{i\ell}
=d_iP_\ell$. Thus, each inventory returns to its initial level at the end of
this period. It remains to choose initial inventories so that
demand is met throughout the period and to bound the resulting
long-run average cost.

\paragraph{Choosing  initial inventories.}
Let $\hat{\pi}_\ell$ be the policy that repeats the modified schedule of period $P_\ell$,
starting with initial inventory $I_i^{\hat{\pi}_\ell}(0)=I_i^\pi(a)+d_i\tau$ for
each commodity $i \in [n]$. To verify feasibility,
let us consider some $t\in[0,P_\ell]$. Let $v$ be the total inserted time elapsed
by $t$, and let $v_i$ be the portion of that time spent producing
commodity $i$. Deleting the inserted intervals identifies $t$
with time $s=a+t-v$ in the original segment. Relative to the
original schedule at $s$, the new schedule starts with $d_i\tau$
additional units, meets $d_i v$ additional units of demand, and
produces $p_i v_i$ additional units. Hence,
\begin{equation} \label{eqn:cyclic-density-proof_1}
I_i^{\hat{\pi}_\ell}(t)-I_i^\pi(s)
 ~~=~~d_i(\tau-v)+p_i v_i
 ~~\in~~[0,p_i\tau],    
\end{equation}
since $0\le v_i\le v\le\tau$ and $d_i<p_i$.
In particular, stockouts cannot occur, and since each
inventory returns to its initial level at time $P_\ell$, repeating
this schedule forms a feasible cyclic policy.

\paragraph{Cost bound.} On the portions retained from $[a,b_\ell)$,
the increase in holding cost rate is at most
$\tau \cdot \sum_{i\in[n]}h_i p_i$, given relation~\eqref{eqn:cyclic-density-proof_1}. To account for the inserted portions,
recall that we insert idle time at the beginning and extend only
the first run of each commodity. Since the endpoint $a$ is fixed, these first
runs all end within a finite interval of the original
schedule, independent of $\ell$. The original inventories are
bounded on this interval, implying
a uniform bound $M$ on the holding cost rate during the inserted
portions. Their total duration is $\tau$, and no additional setup
is required. Consequently,
\[
 {\cal C}(\hat{\pi}_\ell)
 ~~\le~~
 \frac{{\cal C}_{[a,b_\ell)}(\pi)
       +(b_\ell-a)\tau \cdot \sum_{i\in[n]}h_i p_i+M\tau}{P_\ell}.
\]
Since $P_\ell=b_\ell-a+\tau$ and $(b_\ell-a)_{\ell \in \bbN}$ tends to infinity, we conclude that $\limsup_{\ell\to\infty}{\cal C}(\hat{\pi}_\ell) \le C+\tau \cdot \sum_{i\in[n]}h_i p_i$. Finally, choosing $\tau$ so that $\tau \cdot \sum_{i\in[n]}h_i p_i<\eps$ yields ${\cal C}(\hat{\pi}_\ell)\le {\cal C}(\pi)+\eps$ for all sufficiently large $\ell$, as desired.

\subsection{Proof of Lemma~\ref{lem:individual}}
\label{app:proof_lem_individual}

Focusing on commodity $i$, let us recall that $\frac{ \dd q_i^\pi(t) }{ \dd t } = \ind\{t\in\mathcal{P}_i\}-\rho_i =a_i^\pi(t)-\rho_i$  almost everywhere. Therefore,
\begin{align}
\int_{\mathcal S} q_i^\pi(t)\left(a_i^\pi(t)-\rho_i\right) \dd t
 ~~=~~ \int_{\mathcal S} q_i^\pi(t) \cdot \frac{ \dd q_i^\pi(t) }{ \dd t } \dd t ~~&=~~ \left. \frac{ (q_i^\pi(t))^2 }{ 2 } \right]_{ t_0 }^{ t_0 + P} \nonumber \\
 ~~&=~~ \frac{(q_i^\pi(t_0+P))^2-(q_i^\pi(t_0))^2}{2}
 ~~=~~ 0,  \label{eqn:lem_individual_1} 
\end{align}
where the last equality holds since the policy $\pi$ is cyclic, with period length $P$. Noting that commodity $i$ is
produced for $\rho_iP$ units of time in a single period, its mean stock across production periods must be equal to its mean stock overall, $Q_i^\pi$, since by equation~\eqref{eqn:lem_individual_1},
\begin{equation} \label{eqn:lem_individual_2}
\frac{1}{\rho_iP} \cdot \int_{\mathcal S}a_i^\pi(t)q_i^\pi(t) \dd t
 ~~=~~ \frac{1}{P} \cdot \int_{\mathcal S}q_i^\pi(t) \dd t
 ~~=~~ Q_i^\pi .  
\end{equation}

Now consider a production run of commodity $i$ occupying $[x,x+\ell)$.
At elapsed time $u\in[0,\ell]$, its normalized inventory is
$q_i^\pi(x+u)=q_i^\pi(x)+(1-\rho_i)u \geq (1-\rho_i)u$. Consequently, 
\begin{equation}
\int_x^{x+\ell}q_i^\pi(t)\dd t ~~=~~ \int_0^\ell q_i^\pi(x+u) \dd u
 ~~\geq~~ \left. (1-\rho_i) \cdot \frac{ u^2}{2}  \right]_0^{ \ell } ~~=~~ \frac{(1-\rho_i)\ell^2}{2}. \label{eqn:lem_individual_3}   
\end{equation}
We conclude the desired lower bound on $Q_i^\pi$ by observing that
\begin{align}
 Q_i^\pi
~~&=~~ \frac{ 1 }{ \rho_iP} \cdot  \int_{\mathcal S}a_i^\pi(t)q_i^\pi(t) \dd t \label{eqn:lem_individual_4} \\
 &\ge~~ \frac{ 1 }{ \rho_iP} \cdot \frac{1-\rho_i}{2} \cdot \sum_{r \in [m_i^\pi]}(\ell_{ir}^\pi)^2 \label{eqn:lem_individual_5} \\
 &\ge~~ \frac{ 1 }{ \rho_iP} \cdot \frac{1-\rho_i}{2m_i^\pi}
       \left(\sum_{r \in [m_i^\pi]}\ell_{ir}^\pi\right)^2 \label{eqn:lem_individual_6} \\
&=~~ \frac{(1-\rho_i) \rho_i}{2f_i^\pi} . \label{eqn:lem_individual_7}
\end{align}
Here, equality~\eqref{eqn:lem_individual_4} is precisely~\eqref{eqn:lem_individual_2}. Inequality~\eqref{eqn:lem_individual_5} follows from~\eqref{eqn:lem_individual_3}. Inequality~\eqref{eqn:lem_individual_6} holds since $\sum_{r \in [m_i^\pi]}(\ell_{ir}^\pi)^2 \geq \frac{ 1 }{ m_i^\pi } \cdot (\sum_{r \in [m_i^\pi]}\ell_{ir}^\pi)^2$, by Cauchy--Schwarz. Finally, equality~\eqref{eqn:lem_individual_7} is obtained by substituting $\sum_{r \in [m_i^\pi]} \ell_{ir}^\pi = \rho_iP$ and $\mu_i^\pi= \frac{\rho_i}{f_i^\pi}$.

\subsection{Proof of Lemma~\ref{lem:block-inventory}} \label{app:proof_lem_block-inventory}

\paragraph{Proof of item~1.} To bound the contribution of all other commodities'
blocks to $Q_j^\pi$, we observe that
\[ \sum_{i\ne j}Y_{ij}^\pi ~~=~~
\frac{1}{P} \cdot \int_{\mathcal S}
\left(\sum_{i\ne j}a_i^{+\pi}(t)\right)q_j^\pi(t)\dd t ~~\le~~ \frac{1}{P} \cdot \int_{\mathcal S}
\left(1- a_j^\pi(t)\right)q_j^\pi(t)\dd t ~~=~~ (1-\rho_j)Q_j^\pi. \]
The inequality above follows by noting that  $\sum_{i\ne j}a_i^{+\pi}(t) + a_j^\pi(t)\le 1$, since at most one block is in progress at any time,
and blocks of other commodities cannot be in progress while $j$ is
being produced. The last equality is obtained by plugging in the definition of $Q_j^\pi$ and equation~\eqref{eqn:lem_individual_2}.

\paragraph{Proof of item~2.} For any pair of commodities $i \neq j$,
\begin{align}
 Y_{ij}^\pi+Y_{ji}^\pi
 ~~&=~~
 \frac{1}{P} \cdot \int_{\mathcal S}
 \left(a_i^{+\pi}(t)q_j^\pi(t)
      +a_j^{+\pi}(t)q_i^\pi(t)\right)\dd t \nonumber \\
 &\ge~~
 \frac{1}{P} \cdot \int_{\mathcal S}
 \left(a_i^\pi(t)q_j^\pi(t)
      +a_j^\pi(t)q_i^\pi(t)\right)\dd t \label{eqn:proof_lem_block-inventory_eq1} \\
 &=~~
 \frac{1}{P} \cdot \int_{\mathcal S}
 \left( \left( \rho_i+\frac{\dd q_i^\pi(t)}{\dd t} \right) q_j^\pi(t)
      +\left( \rho_j+\frac{\dd q_j^\pi(t)}{\dd t} \right)q_i^\pi(t)\right)\dd t \label{eqn:proof_lem_block-inventory_eq2} \\
&=~~
 \rho_i Q_j^\pi+\rho_j Q_i^\pi
 +\frac{1}{P}  \cdot\int_{\mathcal S}
 \frac{\dd}{\dd t}\left(q_i^\pi(t)q_j^\pi(t)\right)\dd t \nonumber \\
 &=~~
 \rho_i Q_j^\pi+\rho_j Q_i^\pi \left. +\frac{q_i^\pi(t)q_j^\pi(t)}{P} \right]_{t_0}^{ t_0+P} \nonumber \\
 &=~~ \rho_i Q_j^\pi+\rho_j Q_i^\pi. \label{eqn:proof_lem_block-inventory_eq3}
\end{align}
Here, inequality~\eqref{eqn:proof_lem_block-inventory_eq1} holds since each production interval is contained
in its corresponding block, and therefore $a_i^{+\pi}(t) \geq a_i^{\pi}(t)$ for every commodity. To obtain equality~\eqref{eqn:proof_lem_block-inventory_eq2}, we substitute
$a_i^\pi(t)=\rho_i+\frac{\dd q_i^\pi(t)}{\dd t}$ and the analogous
identity for $j$, both valid almost everywhere. Finally, equality~\eqref{eqn:proof_lem_block-inventory_eq3} follows by noting that $q_i^\pi(t_0) = q_i^\pi(t_0+P)$ and $q_j^\pi(t_0) = q_j^\pi(t_0+P)$, since the policy $\pi$ has period length $P$.

\paragraph{Proof of item~3.}
Fixing a pair of commodities $i \neq j$, let us consider an $i$-block
occupying $[v,v+x)$. Since $j$ is not produced during this block,
its stock at any point must cover at least the demand for the
remainder of the block, and therefore $q_j^\pi(v+u) \ge \rho_j(x-u)$ for every $u\in[0,x]$. By integrating over the entire block, we have
\begin{equation} \label{eqn:proof_lem_block-inventory_eq4}
\int_v^{ v+x } q_j^\pi(t) \dd t ~~\ge~~ \rho_j \cdot \int_0^x (x-u)\dd u  ~~=~~ \frac{\rho_j x^2}{2}.    
\end{equation}
Given this relation, we obtain a lower bound on $Y_{ij}^\pi$ by observing that 
\begin{align}
Y_{ij}^\pi  ~~&=~~ \frac{1}{P} \cdot \int_{\mathcal S}
 a_i^{+\pi}(t)q_j^\pi(t)\dd t \nonumber \\
 &\ge~~
 \frac{\rho_j}{2P}  \cdot 
 \sum_{r \in [m_i^\pi]}(s_i+\ell_{ir}^\pi)^2 \label{eqn:proof_lem_block-inventory_eq5} \\
 &\ge~~
 \frac{\rho_j}{2P m_i^\pi} \cdot 
 \left(\sum_{r \in [m_i^\pi]}(s_i+\ell_{ir}^\pi)\right)^2 \label{eqn:proof_lem_block-inventory_eq6} \\
 &=~~
 \frac{\rho_j f_i^\pi}{2}(s_i+\mu_i^\pi)^2 \label{eqn:proof_lem_block-inventory_eq7} \\
 &=~~
 \frac{\rho_j}{2}
 \left(\frac{\rho_i^2}{f_i^\pi}+2\rho_i s_i+f_i^\pi s_i^2\right). \label{eqn:proof_lem_block-inventory_eq8}
\end{align}
Here, inequality~\eqref{eqn:proof_lem_block-inventory_eq5} applies the single-block bound~\eqref{eqn:proof_lem_block-inventory_eq4} with
$x=s_i+\ell_{ir}^\pi$. Inequality~\eqref{eqn:proof_lem_block-inventory_eq6} holds since $\sum_{r \in [m_i^\pi]}(s_i+\ell_{ir}^\pi)^2 \geq \frac{ 1 }{ m_i^\pi } \cdot (\sum_{r \in [m_i^\pi]}(s_i+\ell_{ir}^\pi))^2$ by Cauchy--Schwarz. Equalities~\eqref{eqn:proof_lem_block-inventory_eq7} and~\eqref{eqn:proof_lem_block-inventory_eq8} follow by noticing that $\frac{ 1 }{ m_i^\pi } \cdot \sum_{r \in [m_i^\pi]}(s_i+\ell_{ir}^\pi) = s_i + \mu_i^\pi$, $f_i^\pi = \frac{ m_i^\pi }{ P }$, and $f_i^\pi\mu_i^\pi=\rho_i$.

\subsection{Proof of Lemma~\ref{lem:convex-lower-bound}} \label{app:proof_lem_convex-lower-bound}

Let us consider any feasible cyclic policy $\pi$ of finite cost. We argue that its
frequencies $f^\pi$ and block contributions $Y^\pi$ form a feasible point of \eqref{eq:inventory-convex} whose
objective value is at most $\mathcal C(\pi)$. Consequently, by taking the infimum over all feasible cyclic policies of finite cost
and applying Lemma~\ref{lem:cyclic-density}, we obtain
$\OPT\eqref{eq:inventory-convex}\le\OPTcyc=\OPTdyn$, as required.

First, consider the two terms defining $Q_j(f^\pi,Y^\pi)$ in equation~\eqref{eq:inventory-charge}. By Lemma~\ref{lem:individual}, we have $\frac{\rho_j(1-\rho_j)}{2f_j^\pi} \le Q_j^\pi$ for every commodity $j \in [n]$. 
Moreover, $\sum_{i\ne j}Y_{ij}^\pi \le(1-\rho_j)Q_j^\pi$, by item~1 of Lemma~\ref{lem:block-inventory}. Thus, both terms
in definition~\eqref{eq:inventory-charge} are at most $Q_j^\pi$, and consequently $Q_j(f^\pi,Y^\pi) \le Q_j^\pi$. Now, to verify that $(f^\pi, Y^\pi)$ is a feasible solution to~\eqref{eq:inventory-convex}, we observe that:
\begin{itemize}
    \item Constraint~\eqref{eq:convex-machine}: By inequality~\eqref{eqn:UB_sum_sf}, we have $\sum_{i\in[n]} s_i f_i^\pi \leq \kappa$.

    \item Constraint~\eqref{eq:convex-pair}: By item~2 of Lemma~\ref{lem:block-inventory}, for every $i \neq j$,
    \[ Y_{ij}^\pi +Y_{ji}^\pi ~~\ge~~ \rho_iQ_j^\pi+\rho_jQ_i^\pi ~~\ge~~ \rho_iQ_j(f^\pi,Y^\pi)+\rho_jQ_i(f^\pi,Y^\pi). \]

    \item Constraint~\eqref{eq:convex-block}: By  item~3 of Lemma~\ref{lem:block-inventory}, $Y_{ij}^\pi\ge\frac{\rho_j}{2}
      (\frac{\rho_i^2}{f_i^\pi}+2\rho_i s_i+f_i^\pi s_i^2)$
      for every $i \neq j\in[n]$.
\end{itemize}
In addition, based on  representation~\eqref{eqn:represent_cost_var} of the long-run average cost ${\cal C}( \pi )$, we conclude that the objective value of $(f^\pi, Y^\pi)$ in~\eqref{eq:inventory-convex} is
\[ \sum_{i\in[n]}(K_i f_i^\pi + c_i Q_i(f^\pi, Y^\pi)) ~~\leq~~ \sum_{i\in[n]}(K_i f_i^\pi + c_i Q_i^\pi) ~~=~~ \mathcal C(\pi). \]

\subsection{Approximating the convex relaxation}
\label{app:convex-solvability}

In what follows, we verify the technical conditions needed to apply Theorem~6.3.7 of
\citet{Lovasz2003SDP} in approximating the convex relaxation
\eqref{eq:inventory-convex} within a multiplicative factor of $1 + \eps$.

\paragraph{Strict feasibility.}
For ease of notation, let us define the pair of parameters
\[
 \psi~~=~~\frac{\kappa}{2(1+\sum_{i\in[n]}s_i)}
 \qquad \text{and} \qquad 
 B~~=~~1+\frac{1}{\psi}
   +\max_{i\in[n]}
      \left\{s_i+\frac{\psi s_i^2}{2\rho_i}\right\}.
\]
Next, we create a candidate solution $(f^0, Y^0)$, where $f_i^0=\psi$ for every $i \in [n]$ and $Y_{ij}^0 = B\rho_i\rho_j$ for every $i\ne j$. These are positive rationals of polynomial binary encoding length. To verify the strict feasibility of $(f^0, Y^0)$, note that:
\begin{itemize}
    \item Constraint~\eqref{eq:convex-machine}: By definition of $\psi$, we have $\sum_{i\in[n]} s_i f_i^0 = \psi \cdot \sum_{i\in[n]}s_i < \frac{ \kappa }{ 2 }$.

    \item Constraint~\eqref{eq:convex-pair}: For every $i \neq j$, we have $Y_{ij}^0 +Y_{ji}^0 = 2B\rho_i\rho_j > \rho_iQ_j(f^0,Y^0)+\rho_jQ_i(f^0,Y^0)$. This inequality holds since, by definition~\eqref{eq:inventory-charge},
    \begin{align*}
    \frac{ Q_j(f^0,Y^0) }{ \rho_j } ~~&=~~  \max\left\{ \frac{1-\rho_j}{2f_j^0}, \frac{1}{\rho_j(1-\rho_j)} \cdot \sum_{i\ne j}Y_{ij}^0 \right\} \\
    &\leq~~ \max\left\{ \frac{1}{ \psi }, \frac{B}{1-\rho_j} \cdot \sum_{i\ne j} \rho_i \right\} ~~=~~ \max\left\{ \frac{1}{ \psi }, B \cdot \frac{R - \rho_j}{1-\rho_j}  \right\} ~~<~~ B.
    \end{align*}

    \item Constraint~\eqref{eq:convex-block}: For every $i \neq j$,  note that $B> \frac{ \rho_i }{ 2\psi } + s_i + \frac{ \psi s_i^2 }{ 2\rho_i }$ by definition, and therefore $Y_{ij}^0 = B\rho_i\rho_j > \frac{\rho_j}{2}  (\frac{\rho_i^2}{f_i^0}+2\rho_i s_i+f_i^0 s_i^2)$, as $f_i^0 = \psi$.
\end{itemize}

\paragraph{A bounded epigraph with an interior ball.} Letting $M$ be twice the objective value of $(f^0,Y^0)$, we introduce an epigraph variable $\eta$ and impose $\sum_{i\in[n]} (K_i f_i+c_iQ_i(f,Y)) \le \eta \le M$. This cap preserves $\OPT\eqref{eq:inventory-convex}$, since $(f^0,Y^0)$ already has objective value $\frac{ M }{ 2 }$. It is easy to verify that every point $(f,Y, \eta)$ in the capped epigraph satisfies
\[ f_i ~~\in~~ \left[ \frac{c_i\rho_i(1-\rho_i)}{2M}, \frac{M}{K_i} \right],
 \qquad Y_{ij} ~~\in~~ \left[ 0,\frac{(1-\rho_j)M}{c_j} \right], 
 \qquad \text{and} \qquad \eta ~~\in~~ [0,M], \]
meaning that the capped epigraph is compact and admits an enclosing ball with polynomially encoded rational center and radius.

We next exhibit a ball contained in the capped epigraph, centered
at the point $(f^0,Y^0, \frac{ 3M }{ 4 })$ that strictly satisfies all defining inequalities, and we designate their smallest slack by $\sigma>0$. To this end, if every coordinate
changes by at most $r\le \frac{ \psi }{ 2 }$, each frequency remains positive, since $f_i^0=\psi$, and in addition, $|\frac{1}{f_i}-\frac{1}{\psi}|
 =\frac{|f_i-\psi|}{f_i\psi}
 \le\frac{2r}{\psi^2}$. Since all constraint functions are formed from affine terms,
reciprocals, and maxima, this bound allows us to compute a rational
constant $A\ge1$ such that each constraint function changes by at
most $Ar$. Both $\sigma$ and $A$ have polynomial binary encoding
lengths. Consequently, the Euclidean ball centered at
$(f^0,Y^0, \frac{ 3M }{ 4})$ with radius $r=\min\{\frac{\psi}{2},\frac{\sigma}{2A}\}$ 
lies in the capped epigraph, as every coordinate changes by at most $r$, and every inequality retains slack at least $\frac{ \sigma }{ 2 }$.
This radius is positive, rational, and has polynomial encoding
length. 

\paragraph{Separation and approximation.} A separation oracle must either certify that a given rational point $(f ,Y, \eta)$
belongs to the capped epigraph or return a hyperplane separating
that point from all feasible points. To design a polynomial-time separation oracle, given a rational query point $z = (f ,Y, \eta)$, we first check that
all frequencies are positive. If $f_i\le0$ for some $i \in [n]$, the valid inequality $f_i\ge\frac{c_i\rho_i(1-\rho_i)}{2M}>0$ separates the query point from the capped epigraph. Otherwise,
we evaluate the constraints of~\eqref{eq:inventory-convex},
including $Y\ge0$, together with the epigraph inequalities
$\sum_{i \in [n]} (K_if_i+c_iQ_i(f,Y))\le\eta\le M$.
There are $O(n^2)$ such inequalities, all of which can be evaluated
exactly using rational arithmetic. If all hold, we certify
feasibility. If an inequality $g(z)\le0$
is violated, we compute a rational subgradient $a$ of $g$ at $z$
and return the separating inequality $g(z)+a^\top(x-z)\le0$.
Indeed, convexity implies that every feasible point $x$ satisfies
$g(z)+a^\top(x-z)\le g(x)\le0$, whereas $z$ violates this inequality
since $g(z)>0$. 

Together with the inner and outer balls established above,
this oracle allows us to apply Theorem~6.3.7 of
\citet{Lovasz2003SDP} to minimize $\eta$ within any prescribed additive
accuracy. To obtain a multiplicative guarantee, we observe that the positive rational  $b=\min\{1,K_1,\frac{c_1\rho_1(1-\rho_1)}{2}\}$ is a lower bound on $\OPT\eqref{eq:inventory-convex}$. Indeed, at any feasible point
with $f_1\ge1$, the setup term $K_1f_1$ is at least $K_1$.
When $f_1<1$, the holding term $c_1Q_1(f,Y)$ is
at least $\frac{ c_1\rho_1(1-\rho_1) }{ 2 }$, by definition~\eqref{eq:inventory-charge}. Therefore, minimizing $\eta$ with additive accuracy $\eps b$ yields a rational feasible solution whose objective value is at most $\OPT\eqref{eq:inventory-convex}+\eps b \le (1+\eps) \cdot \OPT\eqref{eq:inventory-convex}$.

\section{Additional Proofs from Section~\ref{sec:construction}} \label{app:proofs_sec_construction}

\subsection{Proof of Lemma~\ref{lem:calendar-fit}} \label{app:proof_lem_calendar-fit}

For $k\ge1$, let $W^{(k)}$ and $t^{(k)}$ respectively denote the total length of all unplaced blocks and the clock value when iteration $k$ of our algorithm begins. We prove by induction that $W^{(k)}\le t^{(k)}$ for all iterations performed. For the base case of $k=1$, just before the first iteration, we have
$t^{(1)}=P$. In addition, 
\[  W^{(1)} ~~=~~\sum_{i\in[n]}N_i e_i ~~=~~P \cdot \sum_{i\in[n]} \frac{ e_i }{ T_i}
 ~~=~~PU~~<~~P ~~=~~t^{(1)}, \]
since $N_i e_i = \frac{ Pe_i }{ T_i } = P u_i$ and $U = \sum_{i \in [n]} u_i < 1$, according to the load bound~\eqref{eq:calendar-load}. For the inductive step, suppose that $W^{(k)}\le t^{(k)}$ when iteration $k$ begins, and that some lot remains
unplaced. We consider the two possible actions in this
iteration:
\begin{itemize}
    \item If our algorithm selects a lot, say of commodity $i$, then $e_i\le W^{(k)}\le t^{(k)}$. Hence, its allocated block $[t^{(k)}-e_i,t^{(k)})$ does not extend to the left of $0$. After placing this block, the remaining work $W^{(k+1)}$ and clock value $t^{(k+1)}$ satisfy
    \[ W^{(k+1)} ~~=~~W^{(k)}-e_i  ~~\le~~ t^{(k)}-e_i ~~=~~t^{(k+1)}. \]

    \item Otherwise, we move the clock backward to the largest remaining due date, denoted by $d$. Since no lot is eligible at $t^{(k)}$, we have $d<t^{(k)}$. No block is placed in this iteration, and every unplaced lot has due date at most $d$. Therefore, 
    \[  W^{(k+1)}  ~~=~~W^{(k)}  ~~\le~~ \sum_{i\in[n]}e_i \cdot \left\lfloor\frac{d}{T_i}\right\rfloor  ~~\le~~ d  ~~=~~t^{(k+1)}, \]
    where the second inequality is obtained by evaluating~\eqref{eq:deadline-work} with $t=d$.
\end{itemize}

The induction shows that every placed block lies within
$[0,P)$. Since the clock moves only  backward and each new
block is placed immediately to its left, blocks cannot
overlap. Moreover, if lot $(i,r)$ is selected in iteration
$k$, its completion time is $C_{ir}=t^{(k)}\le d_{ir}$
by our eligibility condition. Finally, each placement
removes one lot, and each idle step makes at least one
remaining lot eligible, meaning that all lots are eventually placed.

\subsection{Proof of Claim~\ref{clm:Dk_first_UB}} \label{app:proof_clm_Dk_first_UB}

\paragraph{Property 1: Right endpoints of outside blocks.} We first argue that $D_k$ is nonpositive at the right endpoint of every block of each commodity $i>k$. To this end, consider such a block, with right endpoint $b$.
When our placement rule selected this block, the clock was
at $b$. For each $j\in[k]$, let $m_j$ be the number of lots of commodity $j$
still unplaced at that step. These are exactly the lots of $j$ that occupy $[0,b)$, since all previously placed blocks start at or after $b$, and all $m_j$ remaining lots are subsequently placed to the left of the selected block. When  $m_j\ge1$, we must have 
$m_jT_j<b$, since otherwise lot $(j,m_j)$
would have been eligible at time $b$ and would have taken precedence
over the selected lot of commodity $i$. Consequently,
\[  D_k(b)~~=~~ \int_0^b \sum_{j\in[k]} a_j^{+{\pi}}(s) \dd s-v_k b ~~=~~ \sum_{j\in[k]}m_je_j-v_kb        ~~\le~~ b \cdot \left( \sum_{j\in[k]}u_j-v_k \right ) ~~=~~ 0, \]
where the inequality above holds since $m_j < \frac{ b }{ T_j } = \frac{ u_jb }{ e_j }$ when $m_j \geq 1$.

\paragraph{Property 2: Idle intervals.} We next argue that $D_k$ is also nonpositive throughout every idle interval. For this purpose, suppose that our rule leaves $[a,b)$ idle by moving its clock from $b$ to the largest remaining due date, $a$. As before, letting $m_j$ be the number of lots of commodity $j$
still unplaced at this step, we similarly have $m_jT_j\le a$ for every $j\in[k]$. These lots are all placed by time $a$, and no previously placed block starts
before $b$. Consequently, for every $t\in[a,b]$,
\[  D_k(t)~~=~~  \int_0^t \sum_{j\in[k]} a_j^{+{\pi}}(s) \dd s-v_k t ~~=~~ \sum_{j\in[k]}m_je_j-v_kt ~~\le~~ v_k(a-t) ~~\le~~ 0. \]
For the corner case of leaving an idle interval after all lots have been placed,
no block has yet begun and $D_k(t)=-v_kt\le0$.

\paragraph{Bounding $\bar{D}_k$: The area-based argument.} The function $D_k$ is continuous and piecewise linear; its slope is $1-v_k$ during  blocks of prefix commodities and $-v_k$ during any other block or idle interval. Across $[0,P]$, this function starts and ends at $0$, since each commodity $j$ occupies exactly $N_j e_j = \frac{ Pe_j }{ T_j } =Pu_j$ time units over the entire period. Let us consider a maximal interval $(x,y)$ on which $D_k$ is positive.
We have $D_k(x)=D_k(y)=0$, even if $x=0$ or $y=P$. Hence, $D_k$
starts this interval by increasing and ends it by decreasing. By property~2, the machine cannot be idle at any time in $(x,y)$. Thus, immediately before $y$, it is processing a block of some commodity $i>k$. Letting $z$ be the start of this block, we know that $x<z<y$, since $D_k$ initially increases. Moreover, only prefix blocks can occupy time in $(x,z)$, as any outside block occupying time there would end in $(x,z]$, at a point where $D_k$ is positive, contrary to property~1. Thus, the graph of $D_k$ over $(x,y)$ forms a triangle: It rises with slope $1-v_k$ up to height $D_k(z)$ and falls with slope $-v_k$ thereafter. The fall lasts $\frac{ D_k(z) }{ v_k }$ time units and lies within a single block of length $e_i$, implying that $\frac{ D_k(z) }{ v_k }\le e_i$. The area in question, which is charged to the $i$-block starting at $z$, is therefore
\[
 \int_x^y D_k(t) \dd t
 ~~=~~\frac{D_k(z)}{2} \cdot \left(\frac{D_k(z)}{1-v_k}+\frac{D_k(z)}{v_k}\right)
 ~~=~~ \frac{(D_k(z))^2}{2v_k(1-v_k)}
 ~~\le~~\frac{v_ke_i^2}{2(1-v_k)}.
\]

Now, distinct positive
intervals cannot be charged to the same block, since its starting point
lies strictly inside the interval charged to it. As a result, due to having $N_i$ blocks of each commodity $i$, by discarding the nonpositive portions of the integral and summing these charges,
\begin{align}
 \bar{D}_k
 ~~&\le~~ \frac{1}{P}\cdot\sum_{i>k} N_i    \cdot\frac{v_ke_i^2}{2(1-v_k)} \nonumber \\
 &=~~ \frac{v_k}{2(1-v_k)}\cdot\sum_{i>k}\frac{e_i^2}{T_i} \label{eqn:bound_hatD_1} \\
 &\le~~ \frac{v_k}{2(1-v_k)}\cdot        \sum_{i>k} {f}_i^2 T_i (s_i+\mu_i)^2 \label{eqn:bound_hatD_2} \\
 &\le~~ \frac{4v_k}{1-v_k}\cdot          \frac12 \cdot \sum_{i>k}{f}_i(s_i+\mu_i)^2 \label{eqn:bound_hatD_3} \\
 &=~~ \frac{4v_k}{1-v_k}\cdot E . \nonumber
\end{align}
Here, equality~\eqref{eqn:bound_hatD_1} follows by substituting $N_i = \frac{ P }{T_i}$. Inequality~\eqref{eqn:bound_hatD_2} holds since $\ell_i = {f}_iT_i \mu_i$ and $f_i T_i \geq 2$ by relation~\eqref{eq:dilation}, implying that $e_i = s_i+ \ell_i \leq {f}_i T_i ( \frac{ s_i }{ 2 } + \mu_i )$. Finally, inequality~\eqref{eqn:bound_hatD_3} is obtained by recalling that  ${f}_iT_i \leq 4$, again by~\eqref{eq:dilation}.

\subsection{Proof of Claim~\ref{clm:prefix-relaxation}} \label{ap:proof_clm_prefix-relaxation}

We begin by writing the weighted sum on the right-hand side of Claim~\ref{clm:prefix-relaxation} as
\[  \sum_{j\in[k]}\frac{u_jQ_j}{\rho_j}
 ~~=~~
 \sum_{j\in[k]}Q_j
 +\sum_{j\in[k]}\frac{(u_j-\rho_j)Q_j}{\rho_j}, \]
noting that in the second term, the coefficients $\{ u_j-\rho_j \}_{j \in [k]}$ are
nonnegative and sum to $v_k-r_k$. Thus, to obtain the claimed
bound, $E \le \frac{ 1-r_k }{ r_k^2 - r_k v_k + v_k } \cdot \sum_{j\in[k]}\frac{u_j Q_j}{\rho_j}$, it suffices to show that $\frac{ {Q}_j }{ \rho_j } \ge E$ for every $j\in[k]$ and that $\sum_{j\in[k]} {Q}_j\ge \frac{r_k}{1-r_k} \cdot E$.

\paragraph{Bounding $\frac{ {Q}_j }{ \rho_j }$ for  $j\in[k]$.} Since $({f}, {Y})$ is a feasible solution to~\eqref{eq:inventory-convex}, by summing constraint~\eqref{eq:convex-block} over all commodities $i > k$,
\begin{equation} \label{eqn:clm_prefix-relaxation_1}
\sum_{i > k } Y_{ij} ~~\ge~~ \frac{\rho_j}{2} \cdot \sum_{i > k } \left(\frac{\rho_i^2}{f_i}+2\rho_i s_i+f_i s_i^2\right)  ~~=~~ \frac{\rho_j}{2} \cdot \sum_{i > k } f_i (\mu_i + s_i)^2 ~~=~~  \rho_j E,
\end{equation}
where the first equality holds since $\mu_i= \frac{ \rho_i }{ {f}_i }$. As such,  to obtain the desired lower bound on $\frac{ {Q}_j }{ \rho_j }$, we observe that by definition~\eqref{eq:inventory-charge},
\[ \frac{ {Q}_j }{ \rho_j } ~~\geq~~ \frac{1}{\rho_j(1-\rho_j)} \cdot \sum_{i\ne j} {Y}_{ij} ~~\geq~~ \frac{1}{\rho_j} \cdot \sum_{i> k} {Y}_{ij} ~~\geq~~ E. \]

\paragraph{Bounding $\sum_{j\in[k]} {Q}_j$.} By exploiting the feasibility of $({f}, {Y})$ once again, we know by  definition~\eqref{eq:inventory-charge} that ${Q}_j \geq \frac{1}{1-\rho_j} \cdot \sum_{i\ne j} {Y}_{ij}$. As such, to arrive at the desired lower bound, we observe that 
\begin{align}
\sum_{j\in[k]}(1-\rho_j){Q}_j
~~&\ge~~ \sum_{1\le i<j\le k}({Y}_{ij}+{Y}_{ji}) +\sum_{j\in[k]}\sum_{i>k}{Y}_{ij} \nonumber \\
&\ge~~ \sum_{1\le i<j\le k}           (\rho_i{Q}_j+\rho_j{Q}_i) +E \cdot \sum_{j\in[k]} \rho_j \nonumber \\
&=~~ \sum_{j\in[k]}(r_k-\rho_j){Q}_j+r_k E. \label{eqn:bound_sum_Q}
\end{align}
To obtain the second inequality, constraint~\eqref{eq:convex-pair} states that ${Y}_{ij} + {Y}_{ji} \ge \rho_i {Q}_j+\rho_j {Q}_i$ for every $i \neq j$, and in addition, $\sum_{i>k}{Y}_{ij} \geq \rho_j E$, according to inequality~\eqref{eqn:clm_prefix-relaxation_1}. By rearranging~\eqref{eqn:bound_sum_Q} and recalling that $r_k < 1$, we indeed have $\sum_{j\in[k]} {Q}_j\ge \frac{r_k}{1-r_k} \cdot E$.

\subsection{Proof of Claim~\ref{clm:prefix-load}} \label{app:proof_clm_prefix-load}

We begin by showing that $v_k \in [r_k, \frac{ 1 + r_k }{ 2 }]$. To this end, note that
\[  0 ~\le~~ v_k-r_k ~~=~~ \sum_{j\in[k]}( u_j - \rho_j) ~~=~~ \sum_{j\in[k]}\frac{s_j}{T_j} ~~\le~~ \frac12 \cdot \sum_{j\in[k]}s_j{f}_j
 ~~\le~~ \frac{\kappa}{2}
 ~~\le~~ \frac{1-r_k}{2}.
\]
Here, the second inequality holds since ${f}_iT_i \geq 2$   by relation~\eqref{eq:dilation}. The third inequality follows from constraint~\eqref{eq:convex-machine}, stating that $\sum_{i\in[n]}s_i {f}_i \le \kappa$. For the final inequality, $r_k\le R=1-\kappa$.

Now, since $0 < v_k < 1$, the claimed inequality can be equivalently written as
\[ \Psi( v_k, r_k ) ~~=~~ \left( \frac{ r_k^2 }{ v_k(1 - r_k)} + 1 \right) \cdot (1-v_k) ~~\geq~~ \frac{ 1 }{ 1 + \sqrt{2} }. \]
For any fixed $r_k \in (0,1)$, both factors defining $\Psi(v_k,r_k)$ are positive and nonincreasing in $v_k$. It therefore suffices to check the largest admissible value, $v_k=\frac{ 1+r_k }{ 2 }$. At this value,
\[  \Psi\left(\frac{1+r_k}{2},r_k\right) -\frac{1}{1+\sqrt{2}}
 ~~=~~ \frac{1+r_k^2}{2(1+r_k)}-\frac{1}{1+\sqrt{2}}
 ~~=~~ \frac{(r_k-(\sqrt{2}-1))^2}{2(1+r_k)}
 ~~\ge~~ 0. \]

\end{document}